\documentclass[11pt]{article}
\usepackage{amsmath,amssymb,amsthm,amscd}
\newtheorem{Th}{Theorem}
\newtheorem{Cor}{Corollary}
\newtheorem{Prop}{Proposition}
\newtheorem{Lm}{Lemma}
\newtheorem{Def}{Definition}

\theoremstyle{definition}

\newtheorem{rem}{Remark}

\newcommand{\Limpl}{\Longrightarrow}
\newcommand{\A}{{\cal A}}

\newcommand{\D}{{\bf D}}

\newcommand{\Z}{\mathbb{Z}}
\newcommand{\R}{\mathbb{R}}
\newcommand{\N}{\mathbb{N}}

\newcommand{\cs}{{\cal S}}
\renewcommand{\P}{{\cal P}}

\newcommand{\M}{{\cal M}}

\newcommand{\F}{{\cal F}}
\newcommand{\iN}{\*\N_{\infty}}

\newcommand{\B}{{\cal B}}

\newcommand{\e}{\varepsilon}
\renewcommand{\l}{\lambda}
\newcommand{\f}{\varphi}
\newcommand{\Mu}{{\rm M}}
\renewcommand{\a}{\alpha}
\renewcommand{\b}{\beta}
\renewcommand{\d}{\delta}
\renewcommand{\D}{\Delta}
\newcommand{\s}{\sigma}

\newcommand{\sqe}{\subseteq}
\renewcommand{\l}{\lambda}

\newcommand{\la}{\langle}

\newcommand{\range}{\mbox{range}}

\newcommand{\dom}{\mbox{dom}}
\renewcommand{\*}{\,^*\!}
\renewcommand{\o}{\,^\circ\,\!\!}

\newcommand{\all}{\forall}
\newcommand{\ex}{\exists}
\newcommand{\Liff}{\Longleftrightarrow}

\newcommand{\st}{\mbox{st}}

\usepackage{latexsym}
\title{Ergodic theorem for a Loeb space and hyperfinite approximations of dynamical systems}
\author{L.Yu. Glebsky, E.I. Gordon, C.W. Henson}
\date{}
\begin{document}
\maketitle

\begin{abstract}
{\footnotesize Although the G.Birkhoff Ergodic Theorem (BET) is trivial for finite
spaces, this does not help in proving it for hyperfinite Loeb
spaces. The proof of the BET for this case, suggested in \cite{Kam}, works,
actually, for arbitrary probability spaces, as it was shown in \cite{KW}.
In this paper we discuss the reason why the usual approach, based on
transfer of some simple facts about arbitrary large finite spaces
on infinite spaces using nonstandard analysis technique, does not work for the BET.  We show that the the BET for hyperfinite
spaces may be interpreted as some qualitative result for very big finite spaces. We introduce the notion
of a hyperfinite approximation of a dynamical system and prove the existence of such an approximation.
The standard versions of the results obtained in terms of sequences of finite dynamical systems are formulated.}
\end{abstract}

\bigskip

\textbf{Introduction}

\bigskip

Many applications of nonstandard analysis (NSA) to the investigation of infinite structures are based on embedding
these structures in appropriate hyperfinite structures~\footnote{See \cite{Alb, Wolf} for the basic notions of
NSA. As a rule a hyperfinite structure used in this paper can be considered as an ultraproduct of finite structures.}. The latter inherit many properties of finite structures due to the transfer principle of NSA. The effectiveness of the described applications of NSA is based on the fact that many properties can be proved in a more straightforward way for finite structures, than for the corresponding infinite ones.

In the framework of standard mathematics the same idea is implemented by an approximation of an infinite structure
by finite ones, e.g. by an embedding of this infinite structure into an appropriate inductive (projective) limit of finite structures. However, in some problems a construction of inductive and projective limits is either
impossible or inappropriate, while the hyperfinite structures and their nonstandard hulls can be effective.
One can find many such examples in stochastic analysis, functional analysis, harmonic analysis and mathematical physics (see e.g. the books \cite{Alb, Gor, Wolf, GKK} and the bibliography therein).

The book \cite{N1} contains the exposition of Probability Theory based not on the Kolmogorov's axiomatics but on the theory of hyperfinite measure spaces considered by the author as a formalization of the notion of a very big finite probability space. In the author's opinion this approach makes the foundation of probability theory more transparent and more natural when compared with the classical approach.

In this paper we discuss the NSA approaches to the Birkhoff Ergodic Theorem (BET). We remind the formulation of this theorem.

\begin{Th} \label{BET} Let $(X,\Sigma, \nu)$ be a probability space and let $T:X\to
X$ be a measure preserving transformation and $f\in L_1(X)$. Set
$$
A_k(f,T,x)=\frac 1k\sum\limits_{i=0}^{k-1}f(T^ix).
$$

Then

\begin{itemize}

\item[1)] there exists a function $\hat f(x)\in L_1(X)$ such
that $A_k(f,T,x)\to\hat f(x)$ as $k\to\infty$ $\nu$-almost everywhere (a.e.);

\item[2)] the function $\hat f$ is $T$-invariant, i.e. $\hat
f(Tx)=\hat f(x)$, for almost all $x\in X$;

\item[3)] $\int_Xfd\nu=\int_X\hat fd\nu$.
\end{itemize}
\end{Th}

We denote a dynamical system in $X$ by a triple $(X,\nu,T)$ whose terms satisfy the conditions of Theorem \ref{BET}.
In this paper two types of probability spaces are considered:

1) a \emph{Loeb space}, \cite{Alb, Wolf} that is a \emph{standard} probability space
constructed from a hyperfinite set endowed with the counting probability measure;

2) a \emph{Lebesgue space} that is a probability space, which is isomorphic modulo measure $0$ to the segment $[0,1]$ with the standard Lebesgue measure, or to a
countable space with an atomic measure, or to the direct sum of the previous two.

It is known that every compact metric space with a Borel probability measure is a Lebesgue space (see e.g. \cite{KSF}).

We call a dynamical system in a Lebesgue (Loeb) space a \emph{Lebesgue (Loeb) dynamical system}.

The known proofs of Theorem \ref{BET} are not very difficult, but
rather tricky (see e.g. \cite{Br}). Since the specialists in dynamical systems are interested mainly (if not to say, only) in dynamical systems in Lebesgue spaces,
it would be worthwhile to simplify (to make more straightforward) the proof of the BET at least for Lebesgue dynamical systems.

Two attempts to apply NSA to the BET in Lebesgue spaces are known \cite{Kam, Kach 1}.

The NSA approach to a proof of the BET for Lebesgue dynamical systems due to Kamae \cite{Kam} consists of two parts.
\begin{enumerate}
\item The proof of the BET for a Loeb dynamical system.
\item The derivation of the BET for a Lebesgue dynamical system from the BET for a Loeb dynamical system.
\end{enumerate}

It is natural to expect that the proof of the BET for a Loeb space is simpler, than its proof in the general case,
since a Loeb space is constructed from a hyperfinite probability space. The BET for hyperfinite probability spaces
follows from the BET for finite probability by the transfer principle, and for the last class the BET is trivial.

However, the triviality of the BET for finite spaces is not conductive to a simple proof of the BET for Loeb spaces. The proof of the BET for a Loeb space, due to Kamae \cite{Kam}, is not easier, than the known proofs of the general case. In fact, it was shown in \cite{KW} that Kamae's proof is applicable for the general case after a slight modification.

In Section 1 we give some explanation why the BET for finite spaces does not help to prove the BET for Loeb space. We suggest an interpretation of the BET for a Loeb space as some statement in Internal Set Theory \cite{N2} about the stabilization of ergodic means in a hyperfinite space (Theorem \ref{NSBET}). Regarding the notion of a hyperfinite set as a formalization of the notion of a very big finite set, we may consider Theorem \ref{NSBET} as a qualitative statement about the behavior of ergodic means of some specific functions on a very big finite probability space.

In the framework of classical mathematics, the notion of a very big probability space can be formalized using sequences
of finite probability spaces, whose cardinalities increase to infinity. As a rule it is rather difficult to obtain,
or even to formulate, results in probability theory in terms of such sequences. So, dealing with very big finite probability spaces, one prefers to use a transition to infinite probability spaces according to Kolmogorov's axiomatic. However, it is often difficult to interpret statements about standard infinite probability spaces in terms of very big finite probability spaces.

NSA suggests some  methods to reformulate results about hyperfinite structures in terms of sequences of finite structures. The standard "sequence" versions of the results about the behavior of ergodic means in a hyperfinite probability space obtained in Section 1, are discussed in the concluding subsection of this section. It is interesting that Theorem \ref{NSBET} does not have any readable and intuitively clear standard "sequence" version, though its interpretation in terms of very big probability spaces is clear and simple. We describe how the statement of Theorem \ref{NSBET} can be seen in computer experiments.

Another attempt to prove the BET for Lebesgue spaces using NSA is contained in the paper \cite{Kach 1} by A.G. Kachurovskii. This paper is written in the framework of E. Nelson's Internal Set Theory (\textbf{IST}) \cite{N2} without any use of Loeb spaces. The author makes an attempt to deduce Theorem \ref{NSBET} from his theorem  about the fluctuations of ergodic means (see Theorem \ref{Kach1} below). This derivation contains a serious gap that is filled here with the help of Loeb spaces. Actually A.G.Kachurovskii can deduce only Proposition \ref{WeakNSBET} from Theorem \ref{Kach1}. This proposition is equivalent to the BET for a Loeb space. However, it is not clear, whether it is sufficient for the proof of the BET for Lebesgue spaces the approach of \cite{Kach 1}. Theorem \ref{Kach1} is very interesting by itself. Its proof (that is much more difficult, than those of the BET) was transformed in the proof of standard theorem about estimates of $\e$-fluctuations of ergodic m!
 eans \cite{Kach 2}. We show that, unlike Theorem \ref{NSBET}, Theorem \ref{Kach1} has a simple standard "sequence" version (Theorem \ref{stKach}).

\bigskip

In Section 2  we discuss the second part of NSA approach to the BET - the derivation of the BET for Lebesgue
dynamical systems from the BET for Loeb dynamical systems. Though the above discussion shows that the BET for
Loeb dynamical systems is not simpler than for the general case, the both parts of NSA approach to the BET
are interesting in their own right. The second part is interesting from the point of view of approximation of
Lebesgue dynamical systems by finite ones. The approach to approximation suggested here differs from the
most popular approach in ergodic theory based on Rokhlin's Theorem (see e.g. \cite{KSF}). Rokhlin's approximations
have many interesting applications to ergodic theory, especially to problems connected with the entropy of dynamical systems.
The disadvantage of Rokhlin's approximations is caused by the difficulty of their construction that makes them near on impossible to be used in computer simulations of
Lebesgue dynamical systems. In this paper we suggest a definition of approximation of a Lebesgue dynamical system by a hyperfinite dynamical
system (Definition \ref{hypapds} below, see also Definition \ref{hypap}) and prove the existence of such approximation for an arbitrary Lebesgue dynamical system
(Theorem \ref{cycle-appr}). This theorem allows us to prove the universality of Loeb dynamical systems, i.e. the existence of measure preserving epimorphism from a
Loeb dynamical system to a Lebesgue dynamical system. The existence of such epimorphism obviously implies the BET for Lebesgue dynamical systems assuming the BET for Loeb systems to have been already established. The universality of Loeb systems was proved in \cite{Kam}. However, the proof presented here is simpler than one of Kamae and allows to construct hyperfinite approximations for concrete dynamical systems.

The definition of a hyperfinite approximation of Lebesgue dynamical systems can be reformulated in terms of an approximation of Lebesgue dynamical systems by finite dynamical systems in many particular cases, but not in the full generality.  However, for the case of a Lebesgue dynamical system, whose measure space is a compact metric space and the measure preserving transformation is continuous, the standard sequence versions of Definitions \ref{hypap}, \ref{hypapds} and Theorem \ref{cycle-appr} are quite simple (see Definition \ref{st-hypap}, Proposition \ref{seqap4} and Theorem \ref{exseqap}).

As examples, we construct in this paper finite approximations for shifts of the unit interval and for the invertible Bernoulli shift and study the behavior
of ergodic means for these approximations. It is interesting that the ergodic means of approximations of irrational shifts of the unit interval stabilize at the average of a function for all infinite $n$. Theorem \ref{un-erg} states that the same is true for all hyperfinite approximations of any uniquely ergodic dynamical system.

The NSA approach to Rokhlin's finite approximations of Lebesgue dynamical systems will be discussed in another paper. Some preliminary results were announced in
\cite{GHL}.

The authors are grateful to Boris Begun, Karel Hrbacek, Peter Loeb and Edgardo Ugalde for helpful discussions of various parts of this paper.

\section{Ergodic Theorem for Hyperfinite Loeb spaces.}

Recall that a number $\a\in\*\R$ is said to be finite if $|\a|<a$ for some $a\in\R$, otherwise $\a$ is said to be infinite. For every finite $\a\in\*\R$ there exists a unique standard real $\o\a$ such that $\a\approx\o\a$. It is called the standard part or a shadow of $\a$. We use also the following notations:
$\a\sim\infty$ for infinite $\a$, $\a\gg 0$ for $\a\not\approx 0$, $\b\gg\a$ for $\frac{\a}{\b}\approx 0$, $\iN$ for $\*N\setminus\N$.

First we modify a trivial proof of the BET for finite spaces for the
case of hyperfinite spaces.

Throughout this article we consider a hyperfinite set $Y$ such that
$|Y|=M$, where $M$ is an infinite number ($M\in\iN$) and an
internal permutation $T:Y\to Y$. We say that $p(x)$ is the $T$-period of an element $x\in Y$
if its orbit $Orb(x)=\{x,Tx, ..., T^kx....\}$ has exactly $p(x)$ distinct elements. Let $G:Y\to\*\R$ be an internal
function. Denote $A_{p(x)}(G,T,x)$ (the average of $G$ along the orbit of $x$) by $Av_x(G)$.

If there exist an element $x\in Y$ that has a $T$-period $M$, then the $T$-period of any element in $Y$ is equal to $M$.
In this case we say that this permutations is \emph{transitive}, since it generates the transitive cyclic group of permutations of $X$.

If $T$ is transitive, then for every internal function $G:Y\to\*\R$ one has
obviously $\all\, x_1,x_2\in Y\ Av_{x_1}(G)=Av_{x_2}(G)=Av(G)$. So,

\begin{equation} \label{1}
Av(G)=\frac 1M\sum\limits_{y\in Y}G(y)=\frac
1M\sum\limits_{j=0}^{M-1}G(T^jx)=A_M(G,T,x),
\end{equation}
for every $x\in Y$.

Denote the maximal value of $G$ along the orbit of $x$ by $\max_xG$.

\begin{Prop} \label{INTBET} For an arbitrary internal
function $G:Y\to\*\R$ for all $x\in Y$ and for all $N\in\iN$ such that
$\frac{p(x)\cdot\max_x |G|}N\approx 0$ (notation: $N\gg p(x)\cdot\max_x
|G|$) one has $A_N(G,T,x)\approx Av_x(G)$.
\end{Prop}

\textbf{Proof}. It is enough to prove this statement for a transitive permutation $T$. In this case $\all\, x\in Y p(x)=M,\ Av_x(G)=Av(M)$ and $N\gg M\cdot\max |G|$.

Let $N=qM+r$, $0\leq r<M$. Then
$$A_N(G,T,x)=\left(\frac 1{qM+r}-\frac
1{qM}\right)\sum\limits_{i=0}^{qM-1}G(T^ix)+\frac
1{qM}\sum\limits_{i=0}^{qM-1}G(T^ix)+ \frac
1{qM+r}\sum\limits_{i=qM}^{qM+r-1}G(T^ix).
$$
Obviously $\frac 1{qM}\sum\limits_{i=0}^{qM-1}G(T^ix)=\frac
1M\sum\limits_{i=0}^{M-1}G(T^ix)=Av(G)$ by (\ref{1}). The first
and the third terms are infinitesimal since $N\gg M\cdot\max |G|$
$\Box$

We list now the necessary definitions and facts concerning Loeb spaces.
We need here only a particular case of a Loeb space, namely the Loeb space constructed from the
hyperfinite set $Y$ endowed with the uniform probability measure.

Define the external finitely additive measure $\mu$ on the algebra $\P^{int}(Y)$ of internal subsets
of $Y$ by the formula
$$
\mu(B)=\o\left(\frac{|B|}M\right).
$$

The $\aleph_1$-saturation and the Caratheodory theorem imply the
possibility to extend $\mu$ on the $\s$-algebra $\s(\P^{int}(Y))$
generated by $\P^{int}(Y)$. The Loeb space with the underlying set
$Y$ is the probability, space $(Y,P_L(Y),\mu_L)$, where $P_L(Y)$
is the completion of $\s(\P^{int}(Y))$ with respect to $\mu$ and
$\mu_L$ is the extension of $\mu$ on $P_L(Y)$. The measure $\mu_L$
is said to be the \emph{Loeb measure} on $Y$. If necessary we use the notation
$\mu_L^Y$.  We need the following property of the Loeb measure that follows immediately from the
$\aleph_1$-saturation.

\begin{Prop} \label{Loeb-mes} For every set $A\in\P_L$ there exists an internal set $B\sqe Y$ such that
$\mu_L(A\Delta B)=0$. \end{Prop}

\begin{Cor}\label{a.-e.} If $A\in\P_L$, then
$$\mu_L(A)=1\Liff \all^{st}\e>0\,\ex\,B\in\P^{int}(Y)\ (B\sqe A\land\mu(B)>1-\e)$$
\end{Cor}

\begin{rem}\label{Nel-a.-e.} The right hand side of the equivalence in Corollary \ref{a.-e.} can be formalized
in the framework of \textbf{IST} , while this is impossible for the left hand side.
In works based on \textbf{IST} (see e.g. \cite{N1, Kach 1}) the right hand side of the above equivalence is used as
a definition for a certain property to hold almost surely. In this article we say that a property holds $\mu_L-a.e.$ to express the same statement.
\end{rem}

For an arbitrary separable metric space $R$ and an external
function $f:Y\to R$ an internal function $F:Y\to\*R$ is said to be
\emph{a lifting} of $f$ if $\mu_L(\{y\in Y \ |\ F(y)\approx
f(y)\})=1.$

\begin{Prop} \label{mes-func} A function $f:Y\to R$ is measurable
iff it has a lifting.
\end{Prop}

Recall that an internal function $F:Y\to\*\R$ is said to be
\emph{$S$-integrable} if for all $K\in\iN$ one has
\begin{equation}\label{S-int-1}
\frac 1M\sum\limits_{\{y\in Y\ |\ |F(y)|>K\}}|F(y)|\approx 0.
\end{equation}

We need the following properties of $S$-integrable functions.

\begin{Prop} \label{S-int}
1) An $S$-integrable function is almost everywhere bounded.

2) An internal function $F:Y\to\*\R$ is $S$-integrable iff
$Av(|F|)$ is bounded and $\frac 1{M}\sum\limits_{y\in
A}|F(y)|\approx 0$ for every internal $A\sqe Y$ such that
$\frac{|A|}M\approx 0$.

3) An external function $f:Y\to\R$ is integrable w.r.t. the Loeb measure $\mu_L$
iff it has an $S$-integrable lifting $F$, in which
case
$$
\int\limits_Yfd\mu_L=\o Av(F).
$$.
\end{Prop}

We address readers to
\cite{Alb,Wolf} for the proofs of Propositions \ref{Loeb-mes}, \ref{mes-func} and \ref{S-int}.

For the case of $S$-integrable functions one can get some
information about behavior of ergodic means for infinite $N\not\gg
M$, using the same simple considerations.

\begin{Lm} \label{0} If $F:Y\to\*\R$ is an internal function, and $Av(|F|)$ is finite, then for almost all $x\in Y$ $Av_x(|F|)$ is finite. \end{Lm}

\textbf{Proof.} Let $U=\{x\in Y\ |\ Av_x(|F|)\sim\infty\}$. Let $I=\{x_1,\dots x_k\}$ be a \emph{selector}, i.e. it intersects each orbit of $T$ by a single point.
Consider an arbitrary internal set $V\subseteq U$ and show that $\frac{|V|}M\approx 0$. This will prove that $\mu_L(U)=0$. Let $B=\{i\leq k\ |\ Orb(x_i)\cap V\neq\emptyset\}$. Consider the set $V'=\bigcup\limits_{i\in B} Orb(x_i)\supseteq V$. Then $|V'|=\sum\limits_{i\in B}p(x_i)$. To complete the proof of the lemma we
prove that $\frac{|V'|}M\approx 0$.

Let $\l=\min\{Av_x(|F|)\ |\ x\in V\}$. Then $\l\sim\infty$. Suppose that $\frac{|V'|}M\gg 0$. Then
$$
\infty\sim\frac{|V'|}M\cdot\l\leq\frac 1M\sum\limits_{i\in B}p(x_i)\cdot Av_{x_i}(|F|)=\frac 1M\sum\limits_{x\in V'}|F(x)|\leq Av(|F|),
$$
which contradicts the conditions of our Lemma. $\Box$

An internal function $F:Y\to\*\R$ is said to be \emph{bounded} if $\max\{|F(x)|\ |\ x\in Y\}$ is finite. The $\aleph_1$-saturation implies that
$F$ is bounded, if and only if all values of $F$ are finite.

\begin{Th} \label{ErgMeanStab} Let $T:Y\to Y$ be an internal permutation. Then the following statements are true.

1) For every $S$-integrable function $F:Y\to\*\R$, for every
standard positive $a\in\R$, and for every infinite numbers $K,L$
such that $\o\left(\frac KM\right)=\o\left(\frac LM\right)=a$ one
has $A_K(F,T,x)\approx A_L(F,T,x)$ for $\mu_L$-almost all $x\in Y$;

2) The statement 1) holds for all $x\in Y$ such that $\frac{p(x)}M\gg 0$ or 
$\max \{|F(y)|\;|\;y\in Orb(x)\} $ is finite. In particular, the statement 1)
holds for all $x\in Y$ if $T$ is a transitive permutation or $F$ is a bounded function.  

3) If $T$ is a transitive permutation, then for any internal function $F$ such that $Av(|F|)$ is finite ($F$ may be not $S$-integrable here)
and for any $K\in\iN$ such that $\frac KM$ is infinite one has
$A_K(F,T,x)\approx Av(F)$ for all $x\in Y$ .
\end{Th}

\textbf{Proof}. Assume $K>L$ and estimate
$|A_K(F,T,x)-A_L(F,T,x)|$. It is easy to see that
$$
|A_K(F,T,x)-A_L(F,T,x)|\leq\left(\frac 1L-\frac
1K\right)\sum\limits_{k=0}^{L-1}|F(T^kx)|+\frac
1K\sum\limits_{k=L}^{K-1}|F(T^kx)|=U+V.
$$

We need to prove that $U\approx 0$ and $V\approx 0$. 

Let $B=\{T^kx\ |\ k=L,\dots,K-1\}$. Then
$\frac{|B|}M=\frac{K-L}M\approx 0$. Thus, $V=\frac
1M\sum\limits_{y\in B}|F(y)|\approx 0$, due to the
$S$-integrability of $F$ (item 2 of Proposition~\ref{S-int}).

One has
$$
U=\left(\frac ML-\frac MK\right)\frac
1M\sum\limits_{k=0}^{L-1}|F(T^kx)|.
$$
Since $\frac ML\approx\frac MK\approx\frac 1a$ and $\frac
1M\sum\limits_{k=0}^{L-1}|F(T^kx)|\leq\frac
1M\sum\limits_{k=0}^{([a]+1)M-1}|F(T^kx)|$, to establish that $U\approx 0$, it is enough to prove
that the right hand side of the last inequality is finite under conditions 1) and 2) of the theorem.
This is obvious for any $x$ such that $F\upharpoonright Orb(x)$ is a bounded function. 
It is enough to discuss only the case of $[a]=0$. In this case 
$$
\frac 1M\sum\limits_{k=0}^{M-1}|F(T^kx)|=d(x)\cdot\frac 1M\sum\limits_{k=0}^{p(x)-1}|F(T^kx)|+\frac 1M\sum\limits_{k=d(x)\cdot p(x)}^{M-1}|F(T^kx)|=S_1+S_2,
$$
where $M=d(x)\cdot p(x)+r(x)$ and $r(x)<p(x)$. 

Since all terms $S_2$ are distinct, one has $S_2\leq Av(|F|)$, which is finite due to the $S$-integrability 
of $F$. 

If $d(x)$ is finite, then $S_1$ is finite by the same reason. So, $U\approx 0$ in this case. This proves statement 2).

To prove statement 1) notice that
$$
|S_1-Av_x(|F|)|=\frac{r(x)}MAv_x(|F|)
$$
So $S_1$ is bounded for $\mu_L$-almost all $x\in Y$ by Lemma \ref{0}. This proves statement 1).

To prove statement 3) set $K=qM+r,\ r<M$ and notice that $q$ is
an infinite number. Then
\begin{equation} \label{EMS1}
\begin{array}{ll}
A_K(F,T,x)=\frac
1{qM+r}\sum\limits_{i=0}^{qM+r-1}F(T^ix)=\left(\frac1{qM+r}-\frac
1{qM}\right)\sum\limits_{i=0}^{qM-1}F(T^ix)+ \\
+\frac 1{qM}\sum\limits_{i=0}^{qM-1}F(T^ix)+\frac
1{qM+r}\sum\limits_{i=qM}^{qM+r-1}F(T^ix)
\end{array}
\end{equation}
Since $T$ is a cycle of length $M$, one has $\frac
1{qM}\sum\limits_{i=0}^{qM-1}F(T^ix)=Av(F)$ and
\begin{equation} \label{EMS2}
A_K(F,T,x)=\left(\frac{qM}{qM+r}-1\right)Av(F)+Av(F)+ \frac
1{q+\frac rM}\frac 1M\sum\limits_{i=qM}^{qM+r-1}F(T^ix)\approx
Av(F),
\end{equation}
since $Av(|F|)$ is bounded and $|\frac 1{q+\frac rM}\frac
1M\sum\limits_{i=qM}^{qM+r-1}F(T^ix)|\leq\frac 1{q+\frac
rM}Av(|F|)$ $\Box$

\begin{rem} \label{ave} Due to (\ref{1}) one has $A_N(F,T,x)\approx Av(F)$ for $\frac
NM\approx 1$ if $T$ is a transitive permutation. \end{rem}

\begin{Prop} \label{S-int-means} For every $S$-integrable function
$F:Y\to\*\R$ and for every $n\in\*\N$ the function
$F_n=A_n(F,T,\cdot)$ is $S$-integrable and for every $n\in\iN$ the
standard function $\o F_n$ is $T$-invariant a.e.
\end{Prop}
\textbf{Proof}. For $B\subseteq Y$ and for every $n\in\*\N$ one has
\begin{equation} \label{sim1}
\frac 1M\sum\limits_{x\in B}|A_n(F,T,x)|\leq
\frac 1M\sum\limits_{x\in B}A_n(|F|,T,x)=\frac
1n\sum\limits_{i=0}^{n-1}\left(\frac 1M\sum\limits_{x\in
B}|F(T^ix)|\right)
\end{equation}
For $B=Y$ each term in parentheses in (\ref{sim1}) is equal to
$Av(|F|)$ since $T^i$ is a permutation for every $i$. Thus $Av(|A_n(F,T,\cdot)|)\leq
Av(|F|)$, which is finite since $F$ is $S$-integrable.
For $\frac{|B|}{M}\approx 0$ every term in parentheses in the
right hand sum in (\ref{sim1}) is infinitesimal. Thus,
the function $(A_n(F,T,\cdot)$ is $S$-integrable for all
$n\in\*\N$ by item 2) of Proposition~\ref{S-int}.

To prove $T$-invariance of $\o F_n$ for $n\in\iN$ it is enough to show that $G_n(x)=F_n\circ
T(x) - F_n(x)\approx 0$. Obviously $G_n(x)=\frac 1n\left(F(T^{n}x)-F(x)\right)$. Since
$F$ is $S$-integrable, it is bounded a.e. by (\ref{S-int-1}). Thus $G_n(x)\approx 0$ a.e. since $n$ is infinite. $\Box$

\begin{rem} Let $B=\{x\ |\ \all\,n\in\iN\, F_n(x)\approx F_n\circ
T(x)\}$. It is not known if $\mu_L(B)=1$ for every $S$-integrable
function $F$? It is obviously true if $F$ is bounded everywhere.
In this case $B=Y$. If $\o(n/M)>0$, then $\all\,x\in Y\
F_n(x)\approx F_n\circ T(x)$. This easily follows from the
$S$-integrability of $F$.
\end{rem}

\begin{Prop} \label{cont-func} 1) In conditions of the previous proposition
let an internal function $F$  be bounded on $Y$. For every $x\in Y$ define
the standard function $f_x: (0,\infty)\to\R$ by the formula
\begin{equation} \label{cont-func-1}
f_x\left(\o\left(\frac nM\right)\right)=\o F_n(x),
\end{equation}
where $\o\left(\frac nM\right)>0$. Then $f_x$ is well defined and continuous on $(0,\infty)$.

2) The function $\f(x,a)=f_x(a)$ is a measurable function on $Y\times (0,\infty)$, where $Y$ is the Loeb space
and $(0,\infty)$ is the measurable space with Lebesgue measure $dx$.

3) If $T$ is a transitive permutation, then the function $\f(x,a)$ satisfies the following equality
$$\all a\in(0,\infty)\ \ \int\limits_{x\in Y}\f(x,a)dx=\o Av(F).$$
\end{Prop}

\textbf{Proof}. Fix an arbitrary $x\in Y$ and an arbitrary interval
$[a,b]\sqe\R$ such that $a>0$. consider the internal set
$Z_{[a,b]}=\{\frac nM\ |\ n\in\*\N\}\cap [a,b]$ and the internal
function $G_{[a,b]}: Z_{[a,b]}\to\R$ defined by the formula
$G_{[a,b]}\left(\frac nM\right)=F_n(x)$. By Theorem \ref{ErgMeanStab}
the function $G_{[a,b]}$ is $S$-continuous. This proves the first statement
of the proposition.

To prove the second statement, one has to consider the function $G_{[a,b]}$ as an internal function on $Z_{[a,b]}\times Y$.

The third statement follows immediately from the obvious equality $Av(F_n)=Av(F)$ $\Box$

\begin{rem} \label{comp1} The result of Proposition \ref{cont-func} can be observed in computer experiments with finite dynamical systems.
The cardinality $M$, of a finite space $Y$ in these experiments may not be very large, say $M\sim 10^4$. To draw the graph of the function $f_x$
for a finite function $F:Y\to\R$ on an interval $(0,1]$ one has to plot points $\left(\frac nM,A_n(F,T,x)\right)$ for all $n\leq M$. The fact that $F$ is an $S$-integrable function, means that $F$ should not be a $\d$-type function. For example, if $\max |F|\ll M$ (we say in this case that $F$ is bounded), then we may assume that $F$ satisfies the condition of $S$-integrability.

In the Examples 3 and 5 below, the explicit formulas for functions $f_x$ are obtained for some concrete hyperfinite dynamical systems.
\end{rem}

The following example shows that the condition of
$S$-integrability of a function $F$ in Theorem \ref{ErgMeanStab}
is essential and that the assumption of finiteness of
$Av(|F|)$ is not enough in the first part of Theorem
\ref{ErgMeanStab}.

\textbf{Example 1}. Let $Y=\{0,1,\dots, M-1\}$, where $M$ is a
infinite even number. Define the transformation $T:Y\to Y$ by the
formula $T(k)=k+1\,(\mod M)$ and the internal function $F:Y\to\*R$
by the formula $F(k)=M\d_{0k}$ for all $k<M$~\footnote{Here, as usual, $\d_{i,j}=1(0)$, if $i=j(i\neq j)$.}. Then it is easy to see that
$$
A_K(F,T,k)=\left\{\begin{array}{ll} \frac MK,\ &k>M-K \\ 0,\
&0<k\leq M-K\end{array}\right.
$$

We see that if $K<L<M$, and $k$ is such that $M-L<k<M-K$, then
$A_K(F,T,k)\not\approx A_L(F,T,k)$ even if $\frac KM\approx\frac
LM$

Here $Av(|F|)=1$. However, $F$ is not $S$-integrable by Proposition \ref{S-int} (2), since $\frac 1M\sum\limits_{k\in\{0\}}|F(k)|=1$.

For the case of $\frac KM\approx 0$ the BET for Loeb spaces implies the following

\begin{Prop}\label{WeakNSBET} In conditions of Theorem \ref{ErgMeanStab} for any $y\in Y$ 
the $\lim\limits_{n\to\infty} A_n(\o F,T,y)$ exists, if and only if one can find a number $K\in\iN$ 
such that
$$\all\,L\in\iN\ (L<K\Limpl A_K(F,T,y)\approx A_L(F,T,y)\approx\lim\limits_{n\to\infty}A_n(\o F,T,y)).$$
\end{Prop}

The proof of this proposition follows immediately from the following easy
\begin{Lm} \label{Lm1} Let $\{a_n\ |\ n\in\iN\}$ be an internal sequence. Then the standard sequence
$\{\o a_n\ |\ n\in\N\}$ converges if and only if there exist $K\in\iN$ such that $a_K$ is a bounded number
and $\all\,L\in\iN\ (L\leq K\Limpl \o a_L=\o a_K)$ $\Box$
\end{Lm}

The following theorem is stronger then Proposition \ref{WeakNSBET}.

\begin{Th} \label{NSBET}
Let $T:Y\to Y$ be an internal permutation of the Loeb space
$(Y,\P_L(Y),\mu_L)$. Then for every $S$-integrable function
$F:Y\to\*\R$ there exists an infinite hypernatural number $N$ such
that $\mu_L$-almost everywhere for all infinite hypernatural
numbers $L<N$ one has $A_L(F,T,x)\approx A_N(F,T,x)$
\end{Th}

\begin{rem} \label{comp2}  According to Theorem \ref{NSBET} there exists an initial segment
$\A\subseteq\iN\cap\{K\ |\ \frac KM\approx 0\}$ such that ergodic means of $F$ stabilize on $\A$ almost surely.
However, the example 2 below shows that in general the ergodic means $A_K(F,T,k)$ may not stabilize almost surely on the whole set
$\{K\ |\ \frac KM\approx 0\}$. In this case the function $f_x(a)$ corresponding to an $S$-integrable function $F$ can not be extended
by continuity to the point $a=0$. On the graph it can appear as a small cloud in the vicinity of zero. Theorem \ref{NSBET} implies that
there can be found a large enough number $L$, however $L\ll M$ such that for every random number $k\in (0,M)$ the set of points
$\left\{\left(\frac nM, A_n(F,T,k)\right)\ |\ n\leq L\right\}$ looks on display like a horizontal line $y=A_L(F,T,k)$. This change of scale on the x-axis can be interpreted as the observation of neighborhood of zero through a strong microscope. For some finite dynamical systems it is enough to take $M\sim 10^5$ and $L\sim 10^3$
to observe the described phenomenon.
\end{rem}

 To prove Theorem \ref{NSBET} first it is necessary to prove

\begin{Th} \label{a.e.-convergence} Let $f_n:Y\to\R,\ n\in\N$ be a sequence
of measurable functions on $Y$, and $F_n:Y\to\*\R,\ n\in\*N$ be an
\textbf{internal} sequence such that $\all\,n\in\N,\ F_n$ is a
lifting of $f_n$. Then $f_n$ converges to a measurable function
$f$ $\mu_L$-almost everywhere if and only if there exists $K\in\iN$ such that  
 $\mu_L$-almost everywhere $\all\, N\in\iN,\ N<K\Longrightarrow F_N(x)\approx
F(x)$, where $F$ is a lifting of $f$.
\end{Th}

The following lemma is an immediate consequence of the
$\aleph_1$-saturation principle.

\begin{Lm} \label{inf} $\all\,\f:\N\to\iN\,\exists\, N\in\iN\,\all
n\in\N\, N<\f(n)$. $\Box$
\end{Lm}

\textbf{Proof} of Theorem \ref{a.e.-convergence}.

($\Longrightarrow$) Let $f_n$converges to $f$ a.e.  By Egoroff's
theorem
$$
\all k\in\N\,\exists B_k\subseteq Y\, (\mu_L(B_k)\geq 1-\frac
1k\land f_n(x)=\o F_n(x)\ \mbox{converges uniformly on}\
B_k).
$$
WLOG we may assume that $B_k$ is internal, $\frac{|B_k|}{|Y|}\geq
1-\frac 1k$, and $\all\, n,k\in\N\,\all x\in B_k\ F_n(x)\approx f_n(x)$
and $F(x)\approx f(x)$. Then
$$
\exists^{st}\f_k:\N\to\N\,\all^{st}
r\,\all^{st}m>\f_k(r)\max\limits_{x\in
B_k}\,|F_m(x)-F(x)|<\frac 1r.
$$

Consider the internal set
$$
C^k_r=\{N\in\*\N\ |\ \all\,m \,(N>m>\f_k(r)\Limpl\all x\in B_k\
|F_m(x)-F(x)|<\frac 1r)\}
$$
The previous statement shows that $C^k_r$ contains all standard
numbers that are greater that $\f_k(r)$. Thus, there exists infinite $N^k_r\in
C^k_r$. By Lemma \ref{inf} $\exists K\in\iN\, \all^{st} k,r\,
K<N^k_r$. Obviously, this $K$ satisfies Theorem \ref{a.e.-convergence}

($\Longleftarrow$) Let $B=\{x\in Y\ |\ \all N\in\iN (N\leq K\Limpl
F_N(x)\approx F(x)\},\ A_n=\{x\in Y\ |\ f_n(x)\approx F_n(x)\},\
n\in\N,\ A=\{x\in Y\ |\ F(x)\approx f(x)\}, \ C=B\cap
A\cap\bigcap\limits_{n\in\N}A_n.$

By conditions of the theorem $\mu_L(C)=1$. Fix an arbitrary $x\in C$, and an arbitrary
$r\in\N$. The internal set $D=\{n\in\*\N\ \ |\ |F_n(x)-F(x)|\leq \frac 1r\}$ contains
all infinite numbers that are less or equal to $K$. So $\exists n_0\in\N\, \all n>n_0\
|F_n(x)-F(x)|\leq \frac 1r\}$. Since $F_n(x)\approx f_n(x)$, the same holds for
$f_n(x)$ and $\o F(x)$. Thus, $f_n$ converges to $f=\o F$ a.e.     $\Box$.

\textbf{Proof of Theorem \ref{NSBET}.} In conditions of Theorem
\ref{NSBET} let $f=\o F$ and $f_n(x)=A_n(f,T,x),\ n\in\N$ and
the internal sequence $F_n(x)=A_n(F,T,x),\ n\in\*\N$. Then $f\in
L_1(\mu_L)$ and $F_n$ is an $S$-integrable lifting of $f_n$ for
all $n\in\N$. By Theorem \ref{BET} $f_n$ converges to an
integrable function $\hat f$ a.e. Let $\hat F$ be an
$S$-integrable lifting of $\hat f$. Then by Theorem
\ref{a.e.-convergence} there exists $K\in\iN$ such that
$\mu_L$-almost surely $\all N\in\iN\ N<K\Longrightarrow
F_N(X)\approx \hat F(x)$ $\Box$.

\begin{rem} Theorem \ref{a.e.-convergence} also makes it possible to deduce
Theorem \ref{BET} for the Loeb space $Y$ and its internal permutation
$T$ from the Theorem \ref{NSBET}. We use the
same notations as in the proof of Theorem \ref{NSBET}. If there
exist $K\in\iN$ such that a.e. $\all N\in\iN\ N<K\Limpl
F_N(x)\approx\F_K$, then by Theorem \ref{a.e.-convergence} the
sequence $f_n$ converges a.e. to $f=\o F_K$. By Proposition
\ref{S-int-means} $f$ is $T$-invariant a.e. and $F_K$ is
$S$-integrable. So $f$ is integrable $\Box$
\end{rem}

The following example shows, that in general the result of Theorem
\ref{NSBET} can not be extended on the whole galaxy $\{K\in\iN\ |\
\frac KM\approx 0\}$.

\textbf{Example 2} Consider the same Loeb space $Y$ and the same
transformation $T$ as in Example 1. Fix an infinite $K\in Y$ such
that $\frac KM\approx 0$. Let $M=RK+S$, where $0\leq S<K$.
Define the internal function $G:Y\to\*\R$ by the formula
$$
G(k)=\left\{\begin{array}{ll} 1,\ &mK\leq k<(m+1)K,\ m<R,\ m\ \mbox{is even} \\
0,\ &(mK\leq k<(m+1)K,\ m<R,\ m\ \mbox{is odd})\lor RK\leq m<M
\end{array}\right.
$$
The internal function $G$ is bounded and, thus $S$-integrable.

Fix an arbitrary $N\in\iN$ such that $N/K\approx 0$ and consider
the set $B=\{k\in Y |\ \forall n\leq N\ G(T^n(k))=G(k)\}$. It is
easy to see that $\mu_L(B)=1$. Indeed, obviously $B\supseteq
A=\bigcup_{1\leq n\leq R}\{(n-1)K,\dots,nK-1-N\}$
and $\mu_L(A)=\o\left(\frac RM (K-N-1)\right)=1$. Now it is easy to
see that $\all L\leq N\,\all k\in B\ A_L(G,T,k)=A_N(G,T,k)=G(k)$.
Thus, $N$ satisfies Theorem \ref{NSBET}.

Let $D=\{k\in Y\ | mK\leq k<mK+\frac K2\}$. Then, $\mu_L(D)=\frac
12$. It is easy to see, that $\all\, k\in D\, \all n\leq\frac K2\
G(T^n(k))=G(k)$, thus, $A_{\left[\frac K2\right]}(G,T,k)=G(k)$.

To show that Theorem \ref{ErgMeanStab} fails for the galaxy
$\{K\in\iN\ |\ \frac KM\approx 0\}$ it is enough to prove that
$A_K(G,T,k)\not\approx A_{\left[\frac K2\right]}(G,T,k)$ almost
everywhere on $D$. For every standard $n$ consider the set
$D_n=\{k\in D\ |\ |A_K(G,T,k)-G(k)|<\frac 1n\}$. It is enough to
prove that $\lim\limits_{n\to\infty}\mu_L(D_n)=0$. Since the
cardinality of the set $D_n\cap [mK,(m+1)K)$ is the same for all
$m<R$, it is enough to calculate the cardinality of $E_n=D_n\cap
\left[0,\frac K2\right]$. Recall that for $k\in E_n$ one has
$A_{\left[\frac K2\right]}(G,T,k)=G(k)=1$. On the other hand
$A_K(G,T,k)=\frac{K-k}K=1-\frac kK$. So, $|E_n|=\left[\frac
Kn\right]$ and $\mu_L(D_n)=\o\left(R\left[\frac
Kn\right]\right)/(RK+S)\to 0$. $\Box$

The following example shows that one actually can obtain
essentially distinct means for distinct $a=\o\left(\frac
KM\right)$.

\textbf{Example 3}

Let $Y,T$ be the same as in the previous examples and $F(k)=\frac
kM$. For $0\leq a\leq 1$ define the (standard) function $f:[0,1]^2\to\R$
by the formula
$$
\psi_a(t)=\left\{\begin{array}{l} t+\frac a2,\ 0\leq t\leq 1-a
\\ t+\frac a2-1+\frac 1a(1-t),\ 1-a<t\leq 1\end{array}\right.
$$
Then
$$
\o A_K(F,T,k)=\left\{\begin{array}{l} \psi_0\left(\o\left(\frac kM\right)\right),\
\frac KM\approx 0,\ \forall k<M-K \\
\ \psi_a\left(\o\left(\frac kM\right)\right)\ \frac KM\approx a>0,\
\forall k\in X
\end{array}\right.
$$

Notice, that in the first case stabilization holds almost surely
in accordance with Theorem \ref{NSBET}, while in the
second one it holds for all $k\in X$ in accordance with the
Theorem \ref{ErgMeanStab}.

The function $\psi_1(x)=\frac 12$ for all $x$. This is the average of
$f(x)=x$ on $[0,1]$. Notice that $\*f\left(\frac kM\right)=F(k)$.

For $a>0$ one has $\psi_a\left(\o\left(\frac kM\right)\right)=\f(k,a)$, where $\f(x,a)$
is the function defined in Proposition \ref{cont-func}.

\bigskip

\centerline{\emph{Remarks on A.G. Kachurovskii's paper \cite{Kach 1}}}

\bigskip

In the paper \cite {Kach 1} A. G. Kachurovskii proves the NSA version of his theorem about fluctuations of ergodic means\footnote{See \cite{Kach 2} for the standard version
of this theorem. In fact, the NSA version was the first one.} and undertakes an attempt to deduce the BET for Lebesgue spaces from this theorem. However, his approach to the BET contains one gap, that we are going to discuss here.

Let $N\in\iN$. Following \cite{N1} we say that an internal sequence $s:\{1,\dots,N\}\to\*\R$ is convergent if there exists $\xi\in\*\R$ such that
$\all\, L\in\iN\ (L\leq N\Limpl s(L)\approx\xi)$.

For example, Theorem \ref{NSBET} states the existence of a number $N\in\iN$ such that the sequence $\{A_n(F,T,x)\}_{n=1}^N$ is convergent for $\mu_L$ almost all $x\in Y$.

We say that a sequence $s$ admits $k$ $\e$-fluctuations if there exists a sequence of indices $n_1<n_2<\dots<n_{2k-1},n_{2k}$ such that
$\all, i\leq k\ |s_{2i-1}-s_{2i}|\geq\e$. The following definition and proposition are due to E. Nelson \cite{N1}.

\begin{Def} \label{Nels1} We say that an internal sequence $s$ is a sequence of limited fluctuation, if for every standard $\e>0$ and for every $K\in\iN$ the sequence
$s$ does not admit $K$  $\e-fluctuations$. \end{Def}

\begin{Prop}[\cite{N1} Theorem 6.1] \label{Nels2} If $s$ is a sequence of limited fluctuation then there exists $L\in\iN$ such that the sequence
$s\upharpoonright\{1,\dots,L\}$ is convergent. \end{Prop}

The following theorem is an obvious reformulation of Theorem 1A of \cite{Kach 1}.

\begin{Th} \label{Kach1} If $(Y,\mu_L,T)$ is a transitive Loeb dynamical system, $F:Y\to\*\R$ is an $S$-integrable function and $N\in\iN$, then
the sequence $\{A_n(F,T,x)\}_{n=1}^N\}$ is a sequence of limited fluctuation $\mu_L$-a.e. \end{Th}

Using Proposition \ref{Nels2} one immediately deduces Proposition \ref{WeakNSBET} from Theorem \ref{Kach1}

In \cite{Kach 1} it is stated without proof that Proposition \ref{WeakNSBET} implies Theorem \ref{NSBET}. However, this implication requires
Theorem \ref{a.e.-convergence}. Indeed, Proposition \ref{WeakNSBET}, implies the $\mu_L$-a.e. convergence of the standard sequence of functions
$\{\o A_n(F,T,x)\ |\ n\in\N\}$, which by Theorem \ref{a.e.-convergence} implies Theorem \ref{NSBET}. The exposition in \cite{Kach 1} is developed following book \cite{N1}, i.e. it is based on Nelson's Internal Set Theory (IST) \cite{N2}. On the other hand the proof of Theorem \ref{a.e.-convergence}
is based on properties of external Loeb-measurable functions on $Y$. This theorem cannot be even formulated in IST.
Is it possible to prove Theorem \ref{NSBET} in IST with (or without) the help of Theorem \ref{Kach1}, which is a theorem of IST?
We suggest that the answer is positive. We are obliged to K. Hrbacek for the following arguments in support of this conjecture. In \cite{KR} the Hrbacek set theory (HST)
for nonstandard analysis was introduced. This theory is a modification of Hrabcek's theory of external sets \cite{Hr}. It is proved in Chapter 6 of \cite{KR} that
the $\s$-algebra $Bor(H)$ generated by the algebra of internal subsets of a hyperfinite set $H$ can be defined in HST. Thus, Theorem \ref{a.e.-convergence} can be
formulated and it is reasonable to suggest that it can be proved HST as well as the deduction of Theorem \ref{NSBET} from Theorem \ref{a.e.-convergence}. It is proved in
\cite{KR} that every theorem about internal sets that can be proved in HST can be proved in IST. We do not complete these arguments here, since in any case,
these facts about HST are very non-trivial, so it is hard to expect the existence of a proof of Theorem \ref{NSBET} as simple as the one based on
application of Loeb measures.

In the previous paragraph we discussed the proof of the implication
\begin{equation} \label{A}
 \mbox{Proposition \ref{WeakNSBET}}\quad \Limpl\quad \mbox{Theorem \ref{NSBET}}.
\end{equation}

Let us look on this implication from the point of view of Hyperfinite Descriptive Set Theory (HDST).
The basic definitions, notations and fact concerning HDST can be found e.g. in Chapter 9 of \cite{KR}.

Let $Y$ and $M$ be as above and $H=\{0,1,...,M\}$. Consider the set
\begin{equation} \label{E}
U=\{(y,N)\in Y\times (H\cap\iN)\ |\ \all\, L\in\iN\ (L\leq N\Limpl A_L(F,T,y)\approx A_N(F,T,y))\}
\end{equation}
It is easy to see that $U$ is a $\Pi^1_1$-set.

Implication (\ref{A}) is equivalent to the implication
\begin{equation} \label{B}
\mu_L\left(\{y\in Y\ |\ \ex\, N\in\iN\ (y,N)\in U\}\right)=1\Limpl\ex\,N\in\iN\ \mu_L\left(\{y\in Y\ |\ (y,N)\in U\}\right)=1.
\end{equation}

The following proposition easily follows from Theorem 1.1 of \cite{HR}
\begin{Prop} \label{C} If an arbitrary $\Sigma^1_1$-set $U\sqe Y\times (H\cap\iN)$ satisfies the following property
\begin{equation}\label{D}
\all\,(y,N)\in U\,\all L\in\iN\ (L<N\Limpl (y,L)\in U),
\end{equation}
then the implication (\ref{B}) holds for $U$.
\end{Prop}

\textbf{Proof}. Let $\pi_Y:Y\times H\to Y$ and $\pi_H:Y\times H\to H$ be canonical projections. Suppose that $U$ is a $\Sigma_1^1$-set that satisfies
the antecedent of the implication (\ref{B}). Then, by Theorem 1.1 \cite{HR}, for every $n\in\N$ there exist an internal set $U_n\sqe U$ such that
$\mu_L(\pi_Y(U_n))>1-\frac 1n$. Let $N_n=\min\pi_H(U_n)$. Then, by Lemma \ref{inf}, there exists a infinite number $N$ that is less, than all $N_n$. Since
$U$ satisfies the property (\ref{D}), this $N$ satisfies the consequent of the implication (\ref{B}). $\Box$

Unfortunately, the set $U$ defined by formula ($\ref{E}$) is a $\Pi_1^1$-set, but not a $\Sigma^1_1$-set.

So, it is natural to ask the following

\textbf{Question.} Does an arbitrary $\Pi_1^1$-set $U$ that has the property (\ref{D}) satisfy the implication (\ref{B})? More generally, under what conditions
a set $U$ that has a property $(\ref{D})$ satisfies the implication (\ref{B})?

A sufficient condition is given by the following theorem, that may be of interest also by itself.

\begin{Th} \label{Hens} Let $X$ be an internal set and $\l$ a finitely additive internal probability measure on $X$. Let $f:X\to\iN$ be a countably determined (CD)
function. Then there exist a set $S\sqe X$ such that $\l_L(S)=1$, where $\l_L$ is the Loeb measure associated with $\l$, and a number $K\in\iN$ such that
$f(x)\geq K$ for all $x\in S$.
\end{Th}

\begin{Lm}\label{Hens1} Let $Y\sqe X$ be an internal set and $g:Y\to\*\N$ be an internal function. Then there exists a $\l_L$-null set $Y'\sqe Y$ and a
number $K\in\iN$  such that
\begin{equation} \label{Hens2}
\all\,y\in Y\setminus Y'\ (g(y)\in\iN\Limpl g(y)\geq K).
\end{equation}
\end{Lm}
\newcommand{\wit}{\widetilde}
\textbf{Proof of Lemma \ref{Hens1}}. Let $\wit Y=g^{-1}(\iN)$. Since $\wit Y$ is $\l_L$-measurable, there exists a sequence of internal sets
$\{Y_n\sqe\wit Y\ |\ n\in\N\}$ such that $\l_L(\wit Y\setminus Y_n)\leq \frac 1n$. Take $Y'=\bigcap_{n\in\N}\wit Y\setminus Y_n$, so $\l_L(Y')=0$.
Since $Y_n\sqe \wit Y$, we know that $\all\,y\in Y_n$ $g(y)\in\iN$. Since $Y_n$ and $g$ are internal there
exists $K_n\in\iN$ such that $K_n=\min g(Y_n)$. By Lemma \ref{inf} there exists $K\in\iN$ such that $\all\,n\in\N\ K\leq K_n$. Thus,
for every $y\in (Y\setminus Y')\cap\wit Y$ one has $y\geq K$, since $y\geq K_n$ for some $n$. $\Box$

\textbf{Proof of Theorem \ref{Hens}}. Since $f$ is a CD-function, by Theorem 9.4.7 (i) of \cite{KR} there exists a sequence of internal functions $\{g_n\ |\ n\in\N\}$
such that $f\sqe\bigcup\limits_{n\in\N}g_n$. Without loss of generality we may assume that each $g_n\sqe X\times\*\N$. For every $n\in N$ let $Y_n=\dom(g_n)$. Applying
Lemma \ref{Hens1} to $(Y_n,g_n)$ one gets a $\l_L$-null set $Y'_n\sqe Y_n$ and a number $K_n\in\iN$ such that
$$
\all\,y\in Y_n\setminus Y'_n\ (g_n(y)\in\iN\Limpl g_n(y)\geq K_n).
$$
Again, by Lemma \ref{inf} one obtains $K\in\iN$ such that $K\leq K_n$ for all $n\in\N$.

Set $S=X\setminus\bigcup\limits_{n\in\N}Y'_n,$ so that $\l_L(S)=1$. Since $S\sqe X=\bigcup\limits_{n\in\N}Y_n$, for any $y\in S$ there exists $n$ such that $f(y)=g_n(y)\geq K_n\geq K$. $\Box$

Recall that, if a set $V\sqe A\times B$, then a set $W\sqe V$ is said to be \emph{a uniformization} of $V$, if $W$ is a function and $\dom(V)=\dom(W)$.
If $\nu$ is a $\s$-additive measure on $A$ and $\dom(V)\Delta \dom(W)$ is a $\nu$-null set, then $W$ is a $\nu$-a.e.
uniformization of $V$. The following corollary is clear.

\begin{Cor} Let $(X,\l)$ be the same as in Theorem \ref{Hens}. Let $U\sqe X\times\iN$ satisfy the following conditions:
\begin{enumerate}
\item $\l_L\left(\{y\in X\ |\ \ex\, N\in\iN\ (y,N)\in U\}\right)=1$,
\item $\all\,(y,N)\in U\,\all L\in\iN\ (L<N\Limpl (y,L)\in U)$.
\end{enumerate}
Then $U$ has a countably determined $\l_L$-a.e. uniformization if and only if $$\ex\,N\in\iN\ \l_L\left(\{y\in X\ |\ (y,N)\in U\}\right)=1.\quad\Box$$
\end{Cor}
Unfortunately, the existence of CD-uniformizations is proved only for sets with $\Sigma^0_1$-cross-sections, and there are known examples of sets
with $\Pi^1_1$-cross-sections without CD-uniformizations. So, this corollary is not applicable to the set $U$ defined by (\ref{E}).

\bigskip

\centerline{\emph{Standard versions of the previous results}}

\bigskip

While in nonstandard analysis we use the notion of an infinite number (hyperfinite set) as a formalization of the notion of a very big number (finite set),
in classical mathematics we use the sequences of numbers (finite sets) diverging to infinity to formalize these notions. For example,
in previous sections we considered a hyperfinite set $Y$ and its internal permutation $T:Y\to Y$. If we want to treat the same problems in the
framework of standard mathematics, we have to consider a sequence $(Y_n, T_n)$ of finite sets $Y_n$ whose cardinalities tend to infinity and their permutations $T_n$.
Similarly, internal functions $F:Y\to\*\R$ correspond to sequences $F_n:Y_n\to\R$ in standard mathematics.

First, we discuss what property of such sequences correspond to the property of an internal function $F$ to be $S$-integrable. The following proposition gives a
reasonable answer to this question.

\begin{Prop} \label{uni-int} Let $Y_n$ be a standard sequence of finite sets, such that $|Y_n|=M_n\to\infty$ as $n\to\infty$. Then for an arbitrary sequence
$F_n:X_n\to\R$ the following statements are equivalent:
\begin{enumerate}
\item For every $K\in\iN$ the function $\* F_K$ is $S$-integrable.
\item
\begin{equation}\label{uni-int-1}
\lim\limits_{n,k\to\infty}\frac 1{M_n}\sum\limits_{\{x\in Y_n |\ |F_n(x)|>k\}}|F_n(x)|=0
\end{equation}
\end{enumerate}
\end{Prop}

The proof can be obtained easily by application of the Nelson 's algorithm \cite{N2} to the statement (\ref{S-int-1}).

A sequence $F_n$ that satisfies the statement (2) of Proposition \ref{uni-int} is said to be \emph{uniformly integrable}.

Proposition \ref{uni-int} leads to establishing the standard version of Theorem \ref{ErgMeanStab}.

\begin{Prop} \label{ErgMeanStabSt} In conditions of Proposition \ref{uni-int} let $T_n:Y_n\to Y_n$ be a sequence of transitive permutations and $F_n:Y_n\to\R$ be a uniformly integrable sequence. Consider two sequences of natural numbers $K_n$ and $L_n$ such that $\frac{K_n}{M_n}$ is bounded, $\liminf\frac{K_n}{M_n}>0$ and $\lim\limits_{n\to\infty}\frac{K_n}{L_n}=1$. Then the following two statements  are true.
\begin{enumerate}
\item For any $\e>0$ one has
\begin{equation}\label{ErgMeanStab 0}
\lim\limits_{n\to\infty}\frac 1{M_n}\cdot\left|\{y\in Y_n\ |\ \left|A_{K_n}(F_n,T_n,y)-A_{L_n}(F_n,T_n,y)\right|\geq\e\}\right|=0
\end{equation}
\item. If $T_n$ is a sequence of transitive permutations or $F_n$ is a sequence of bounded functions, then 
\begin{equation}\label{ErgMeanStab1}
\lim\limits_{n\to\infty}\max\limits_{y\in Y_n}\left|A_{K_n}(F_n,T_n,y)-A_{L_n}(F_n,T_n,y)\right|=0
\end{equation}
\end{enumerate}
\end{Prop}

Though the proof of Theorem \ref{ErgMeanStab} can be easily rewritten in (standard) terms of Proposition \ref{ErgMeanStabSt}, we deduce here this proposition from Theorem \ref{ErgMeanStab}, keeping in mind the following much more difficult analysis of the standard version of Theorem \ref{NSBET}. By the same reason we do not use the Nelson's algorithm in our proof. 

\textbf{Proof}. We restrict ourselves to the case of transitive permutations $T_n$. Other cases can be treated in a similar way.

Let $\F$ be a non-principal ultrafilter over $\N$. Consider the nonstandard universe $\*U$ that is the ultrapower of the standard universe $U$ by
$\F$ ($\*U=U^{\N}/\F)$. Let $Y={Y_n}^{\F},\ T={T_n}^{\F},\ F={F_n}^{\F},\
M={M_n}^{\F},\ K={K_n}^{\F},\ L={L_n}^{\F}$ be the classes of corresponding sequences in $\*U$. Then $Y$ is a hyperfinite set of infinite cardinality,
since $\lim\limits_{n\to\infty}M_n=\infty$, thus $\lim\limits_{\F}M_n=\infty$,
$T:Y\to Y$ is an internal transitive permutation on $Y$, $F:Y\to\*\R$ is $S$-integrable by Proposition \ref{uni-int},
$\o\left(\frac KM\right)=\lim\limits_{\F}\frac{K_n}{M_n}>0$, $\o\left(\frac KL\right)=1$, so $\o\left(\frac KM\right)=\o\left(\frac LM\right)>0$.
Thus, the introduced nonstandard objects satisfy all conditions of Theorem \ref{ErgMeanStab} (2). So, one has
$\all\, y\in Y\ |A_K(F,T,y)-A_L(F,T,y)|\approx 0$ and $\max\limits_{y\in Y}|A_K(F,T,y)-A_L(F,T,y)|\approx 0$. Thus,
$$\o\left(\max\limits_{y\in Y}|A_K(F,T,y)-A_L(F,T,y)|\right)= \lim\limits_{\F}\max\limits_{y\in Y_n}\left|A_{K_n}(F_n,T_n,y)-A_{L_n}(F_n,T_n,y)\right|=0.$$
Since, this is true for an arbitrary non-principal ultrafilter $\F$ the equality (\ref{ErgMeanStab1}) is proved. $\Box$

Let us discuss a standard sequence version of Theorem \ref{NSBET}. Assuming that $\{F_n:Y_n\to\R\}$ is uniformly integrable sequence, $|Y_n|\to\infty$ and $T_n$ is a
permutation on $Y_n$, we write in IST the statement that for every infinite $N$ the triple $(Y_N, T_N,F_N)$ satisfies Theorem \ref{NSBET}:
\begin{equation} \label{fnsbet}
\begin{array}{c}\all N\in\iN\,\all^{st}\e>0\,\ex A\sqe Y_N\,\ex L\in\iN\ \left (\frac{|A|}{|Y_N|}\geq 1-\e\land\all K\in\iN\ (K\leq L\Limpl\right.
\\ \left. \Limpl \all y\in A \ A_K(F_N,T_N,y)\approx A_L(F_N,T_N,y))\right)\end{array}
\end{equation}

We call the statement \ref{fnsbet} the IST-sequence version of Theorem \ref{NSBET}

Since only the $\aleph_1$-saturation is needed in the proof of Theorem \ref{NSBET}, to obtain the standard sequence version of this theorem one needs to write
the statement about the validity of (\ref{fnsbet}) in the ultrapower $U^{\F}$ for an arbitrary non-principle ultrafilter $\F\sqe\P(\N)$. We skip this simple exercise.
It is clear, however that the existence of infinite $L$ satisfying certain conditions, means the existence for every non-principal ultrafilter $\F$
\emph{its own} sequence $L_n$ such that $\lim\limits_{\F}L_n=\infty$. So, in this case, unlike the Proposition \ref{ErgMeanStabSt}, it is impossible to get rid of
limits along ultrafilters. This makes the standard version of Theorem \ref{NSBET} very unclear intuitively.

The formalization of Proposition \ref{WeakNSBET} in the framework of IST is not simpler than the sentence (\ref{fnsbet}). To obtain this formalization one should only
place the quantifier $\all y\in A$ in (\ref{fnsbet}) between the quantifiers $\ex A\sqe Y_N$ and $\ex L\in\iN$.

The IST-sequence version of Theorem \ref{Kach1} is much simpler, than those of Theorem \ref{NSBET}. In sake of shortness let us restrict ourselves by sequences of
ergodic means of length $|Y_N|$. It means that we consider only sequences $\{A(F,T,y)_{n=1}^{|Y|}$ in Theorem \ref{Kach1}. Using Theorem \ref{ErgMeanStab}
it is not difficult to deduce the general case of Theorem \ref{Kach1} from this particular case.
Let $Y_n,T_n$ and $F_n$ be as above.
Denote the number of $\e$-fluctuations in the sequence $\{A_n(F_N,T_N,y)\}_{n=1}^{|Y_N|}$ by $Fl(\e,N,y)$.
Let $U(k,\e,N)=\{y\in Y_N\ |\ Fl(\e,N,y)\leq k\}$. Notice, that $U(k,\e,N)$ is an internal set.

Then it is easy that Theorem \ref{Kach1} is equivalent to the following IST-sentence
$$
\all N\in\iN\,\all^{st}\e>0\,\all^{st}\d>0\,\ex^{st}k\ \frac 1{|Y_N|}|U(k,\e,N)|\geq 1-\d
$$
Using Nelson's algorithm it is easy to obtain the following standard sequence version of Theorem \ref{Kach1}
\begin{Th}\label{stKach} Let $Y_n$ be a standard sequence of finite sets, such that $|Y_n|=M_n\to\infty$ as $n\to\infty$,
let $T_n:Y_n\to Y_n$ be a sequence of transitive permutations and let $F_n:Y_n\to\R$ be a uniformly integrable sequence. Then

$$\all\, \e,\d>0\,\ex N,k\,\all\, n>N\ \frac 1{M_n}|U(k,\e,n)|\geq 1-\d.\ \Box$$ \end{Th}

It would be interesting to obtain estimates concerning fluctuations of ergodic means for this theorem similar to those obtained in \cite{Kach 2}.

\bigskip

\section{Hyperfinite approximations of dynamical systems on compact metric spaces}

Let $X$ be a compact metric space with a metric $\rho$. It is known that for every $\xi\in\* X$
there exists a unique standard element $x\in X$ such that $x\approx\xi$. Thus, the external map
$st:\*X\to X$ such that $\all\,\xi\in\*X\ st(\xi)\approx\xi$ is defined. The map $st$ (we also use the notation $st_X$, if necessary) is called
the standard part map. For a hyperreal number $\xi\in\*\R$ obviously one has $\o\xi=\st(\xi)$. The map $st$ is defined in this case not for all
hyperreal numbers but only for bounded ones, since $\R$ is not a compact space.

\begin{Def} \label{hypap} Let $X$ be a compact metric space, $\nu$ be a Borel measure on
$X$ and $\f:Y\to \* X$ be an internal map such that $st\circ\f:Y\to
X$ is a measure preserving map. We say that in this case $(Y,\f)$
is a \emph{hyperfinite approximation (h.a.)} of the measure space $(X,\nu)$.
\end{Def}

In case of $Y\sqe\*X$ and the identical embedding $Y$ we say that $Y$ is a h.a.
of $(X,\nu)$. Obviously, any h.a. $(Y,\f)$ is equivalent to the h.a. $\f(Y)$.

If $R$ is a complete separable metric space and $f:X\to R$ is a
Borel measurable function, then $f\circ st\circ\f:Y\to R$ is a $\mu_L$-measurable
function. A lifting $F:\*X\to\*R$ of the function $f\circ st\circ\f$ is said
to be a lifting of the function $f$. By Proposition \ref{mes-func}
every Borel measurable function $f:X\to R$ has a lifting.

To formulate the standard version of this definition introduce the following notation.
Let $Z\sqe X$ be a finite subset of $X$ and $\d_Z=\frac 1{|Z|}\sum\limits_{z\in Z}\d_z$, where
$\d_Z$ is a Dirac measure at a point $z\in Z$, i.e. $\d_Z$ is a Borel probability measure such that
 for any Borel set $A\sqe X$ one has $\d_z(A)=1\Liff z\in A$.

\begin{Def} \label{st-hypap} In conditions of Definition \ref{hypap} let $\{Y_n\ |\ n\in\N\}$ be a sequence of
finite subsets of $X$. We say that the sequence $Y_n$ approximates the measure space $(X,\nu)$ if the sequence of measures
$\d_{Y_n}$ converges to the measure $\nu$ in the *-weak topology on the space $\M(X)$ of all Borel measures on $X$.
\end{Def}

\begin{Prop} \label{ex-st-hypap}  In conditions of Definition \ref{st-hypap} suppose that every open ball in $X$ has the positive measure $\nu$ and
every set of the positive measure $\nu$ is infinite. Then for every set $A\sqe X$ with $\nu(A)=1$ there exists a sequence $Y_n$ of finite subsets
of $X$ approximating the measure space $(X,\nu)$ such that $\all\,n\in\N\ Y_n\sqe A$.
\end{Prop}

\textbf{Proof}. Since $X$ is a compact metric space, the space $C(X)$ of all continuous functions on $X$ is separable.
Then the space $\M(X)$ of all Borel measures on $X$ is separable in the *-weak topology. By the Krein - Milman Theorem the convex set of probability measures in $\M(X)$
is the closure of the convex combinations of its extreme points, which are Dirac measures. So, the set of all convex combinations with rational coefficients is dense
in the set of probability measures. It is enough to show that for any finite set $E=\{z_1,\dots,z_k\}\sqe X$, for any natural numbers $\{n_1,\dots, n_k\}$,
for any finite set of continuous functions $\{f_1,\dots,f_m\}\sqe C(X)$ and for any $\e>0$, there exists a set $Y=\{y_1,\dots,y_n\}\sqe X$, where $n=n_1+...+n_k$, such that
\begin{equation} \label{inequal}
\all\,i\leq m\ \left|\sum\limits_{j=1}^k\frac{n_j}nf_i(z_j)-\int\limits_Xf_id\d_Y\right|\leq\e.
\end{equation}

Since functions $f_i$ are uniformly continuous on $X$ there exists $\s>0$ such that
$$
\all\/i\leq m\ \all\,x_1,x_2\in X\ \rho(x_1,x_2)<\s\Limpl |f_i(x_1)-f_i(x_2)|<\e
$$
and
$$
\s<\min\{\frac 12\rho(u,v)\ |\ u,v\in E,\ u\neq v\}.
$$
Then $B_{\s}(z_i)\cap B_{\s}(z_j)=\emptyset$ for $i\neq j$. By the conditions of our proposition each set
$B_{\s}(z_i)\cap A$ is infinite. So, there exist a set
$$
Y=\{y_1,\dots,y_{n_1}, y_{n_1+1},\dots,y_{n_2},\dots,y_{n_{k-1}+1},\dots,y_{n_k}\}\sqe A
$$
such that $y_1,\dots,y_{n_1}\in B_{\s}(z_1);\ y_{n_1+1},\dots,y_{n_2}\in B_{\s}(z_2);\ \dots;y_{n_{k-1}+1},\dots,y_{n_k}\in B_{\s}(z_k)$.
Obvious calculation shows that $Y$ satisfies (\ref{inequal}) $\Box$.

\begin{Prop} \label{equivalence} A sequence $Y_n\sqe X$ approximates a measure space $(X,\nu)$ in the sense of Definition \ref{st-hypap} if and only if
for any $N\in\iN$ the set $Y_N$ is a hyperfinite approximation of the measure space $(X,\nu)$. \end{Prop}

\textbf{Proof.} $\Limpl$ Let $Y_n$ approximates $(X,\nu)$ and $N\in\iN$. Then for any $f\in C(X)$ one has
\begin{equation} \label{equivalence1}
\int f(\st(x))d\mu_L=
\o\left(\frac 1{|Y_N|}\sum\limits_{y\in Y_N}\*f(y)\right)=\int\limits_Xfd\nu
\end{equation}
The first equality is due to $\*f$ is a lifting of $f\circ\st$.  
The second follows from Definition \ref{st-hypap} and from the nonstandard analysis definition of the limit of a sequence.
Now $\st:Y_N\to X$ defines a measure $\nu'$ on $X$ that is the image of the Loeb measure of $Y$. Due to (\ref{equivalence1}) and the Riesz 
representation theorem $\nu'=\nu$.

$\Longleftarrow$ Assume that $\st\upharpoonright Y_N:Y_N\to\R$ is a measure preserving transformation for every $N\in\iN$. It is easy to see that for every function
$f\in C(X)$ the internal function $\*f\upharpoonright Y_N$ is a lifting of $f$. So,
$$
\o\left(\frac 1{|Y_N|}\sum\limits_{y\in Y_N}\*f(y)\right)=\int\limits_{Y_N}\o(\*f)d\mu_L=\int\limits_{Y_N}f\circ st\,d\mu_L=\int\limits_Xfd\nu.
$$
Thus, the equality (\ref{equivalence1}) holds for every $N\in\iN$ and by the nonstandard analysis definition of a limit one has
$\lim\limits_{n\to\infty}\int\limits_Xfd\d_{Y_n}=\int\limits_Xfd\nu$ $\Box$

Propositions \ref{ex-st-hypap} and \ref{equivalence} imply the following

\begin{Prop} \label{Space-h-a} Let $\nu$ be a non-atomic Borel measure
on a compact metric space $X$ such that the measure of every ball
is positive. Then for every set $A\sqe X$ such that $\nu(A)=1$
there exists a hyperfinite set $Y\sqe \*A$ such that $(Y,\st)$ is a
h.a. of $(X,\nu)$.
\end{Prop}

As it was mentioned in Introduction, in conditions of Proposition \ref{Space-h-a} the measure space $(X,\nu)$
is isomorphic modulo measure $0$ to the measure space $([0,1],dx)$, where $dx$ is the standard Lebesgue measure.
This means that there exist a set $B\sqe X$ a set $C\sqe [0,1]$ and a bijective map $\psi:B\to C$
such that $dx(C)=\nu(B)=1$ and the maps $\psi,\psi^{-1}$ are measure preserving.

\begin{Lm} \label{hyp-iso} In conditions of the previous paragraph let $Y$ be a h.a. of $(X,\nu))$. Then for every set $D\sqe [0,1]$
with $dx(D)=1$ there exists a bijective lifting $G:Y\to\*[0,1]$ of the map $\psi$ such that
\begin{enumerate}
\item $Z=G(Y)\sqe \*D$;
\item $Z$ is a h.a. of $([0,1],dx)$.
\item $G^{-1}:Z\to \*X$ is a lifting of $\psi^{-1}$.
\end{enumerate}
\end{Lm}

\textbf{Proof}. Let $F:Y\to\*[0,1]$ be a lifting of $\psi$. Let $\s=\frac12\min\{\rho(u,v)\ |\ u,v\in F(Y),\,u\neq v\}$.
Then $0<\s\approx 0$ and $\all\, u\in F(Y)\ B_{\s}(u)\cap F(Y)=\{u\}$. Since $\*\nu (B_{\s}(u))>0$ and $dx(D)=1$ the set $B_{\s}(u)\cap\*D$ contains infinitely
many points, and thus, there exists an internal set $E_u\sqe B_{\s}(u)\cap\*D$ such that $|E_u|=|F^{-1}(u)|$. Establishing bijection between $F^{-1}(u)$ and $E_u$ for every
$u\in F(Y)$, we obtain the bijection $G:Y\to Z\sqe\*D$ that is a lifting of $\psi$. Notice that since $G$ and $G^{-1}$ are bijections they are measure preserving
maps between measure spaces $(Y,\mu_L^Y)$ and $(Z,\mu_L^Z)$.

To prove the second property of the set $Z$, one needs to show that $st_{[01]}\upharpoonright Z:Z\to [0,1]$ is a measure preserving map, i.e. that for every measurable set
$A\sqe[0,1]$ one has
\begin{equation} \label{hyp-iso1}
\mu_L^Z(st^{-1}_{[0,1]}(A)\cap Z)=dx(A)
\end{equation}
One has
\begin{equation}\label{hyp-iso2}
\mu_L^Z(st^{-1}_{[0,1]}(A)\cap Z)=\mu_L^Y(\{y\in Y\ |\ G(y)\in st^{-1}_{[0,1]}(A)\})=\mu_L^Y(\{y\in Y\ |\ \o G(y)\in A\}).
\end{equation}
Since $G$ is a lifting of $\psi$ on has $\o G(y)=\psi(st_X(y))$ for $\mu_L^Y$-almost all $y$. Thus,
\begin{equation}\label{hyp-iso3}
\mu_L^Y(\{y\in Y\ |\ \o G(y)\in A\})=\mu_L^Y(\{y\in Y\ |\ \psi(st_X(y))\in A\})=\nu(\psi^{-1}(A))=dx(A),
\end{equation}
since $st_X\upharpoonright Y:Y\to X$ and $\psi:B\to C\sqe [0,1]$ are measure preserving maps.
The equality (\ref{hyp-iso1}) follows from the equalities (\ref{hyp-iso2}) and (\ref{hyp-iso3}).

To prove the third property of the set $Z$ it is enough to show that $st_Y(G^{-1}(z))=\psi^{-1}(st_{[0,1]}(z))$ for $\mu_L^Z$-almost all $z\in Z$. Since $\psi$ is a bijection,
the last equality is equivalent to the equality $\psi (st_Y(G^{-1}(z)))=st_{[0,1]}(z)$, which follows from the following sequence of equalities  that hold for $\mu_L^Z$-almost
all $z\in Z$:
$$\psi(st_Y(G^{-1}(z))=st_{[0,1]}(G(G^{-1}(z)))=st_{[0,1]}(z). \Box$$

\begin{Def} \label{hypapds} Let $X$ be a compact metric space, $\nu$ be a Borel measure on $X$, $\tau:X\to X$
be a measure preserving transformation of $X$ and $(Y,\f)$ be a h.a. of $(X,\nu)$. Then we say that an internal permutation
$T:X\to X$ is a h.a. of the transformation $\tau$ if the following diagram
\begin{equation} \label{diag}
\begin{CD} Y @>T>> Y \\  @V\f VV  @V\f VV\\ \*X @>\*\tau>>\* X\end{CD}
\end{equation}
is commutative $\mu_L$-a.e.
We say also that the dynamical system $(Y,\mu_L, T)$ is a h.a. of the dynamical system $(X,\nu,\tau)$.

In case of $Y\sqe \*X$ and and an identical imbedding $\f$ the diagram (\ref{diag}) means that $T$ is a lifting of $\tau$ i.e. $\o T(y)\approx\tau(\o y)$ $\mu_L$-a.e.
\end{Def}

\begin{Th} \label{cycle-appr} For every dynamical system
$(X,\nu,\tau)$ such that the measure space $(X,\nu)$ satisfies the conditions of
Proposition \ref{Space-h-a}, and for every h.a. $Y$ of $(X,\nu)$ there exists a h.a. $(Y,\mu_L, T)$.
Moreover, one can choose a h.a. $T$ of $\tau$ to be a transitive permutation.
\end{Th}

\textbf{Proof.} I. Here we prove the existence of a h.a. $(Y,\mu_L,T)$ of the
dynamical system $([0,1],dx,\tau)$.
Let $Y$ be an arbitrary h.a. of the measure space $([0,1],dx)$.  Let $F:Y\to\*[0,1]$ be a lifting of $\tau$.
First we prove  the following statement.

(A) \emph{For every standard $\d>0$ there exists a permutation
$T_{\d}:Y\to Y$ such that}
$$
\frac{|\{y\in Y\ |\ |F(y)-T_{\d}(y)|<\d\}|}M\approx 1.
$$
We deduce (A) from the Marriage Lemma.
Fix a standard $\d>0$ and for every $y\in Y$ set
$S(y)=\*(F(y)-\d,F(y)+\d)\cap Y$. Let $I$ be an arbitrary internal
subset of $Y$. Set $S(I)=\bigcup\limits_{y\in I}S(y)$ and
$B(I)=\bigcup\limits_{y\in I}\*(F(y)-\d,F(y)+\d)$. So,
$S(I)=B(I)\cap Y$. The internal set $B(I)$ can be represented as a
union of a hyperfinite family of disjoint intervals. Since the
length of each of these intervals is not less than $2\d$, their
number is actually finite. Let
$B(I)=\bigcup\limits_{i=1}^n(\xi_i,\eta_i)$, where intervals
$(\xi_i,\eta_i)$ are pairwise disjoint and $n$ is standard.

Consider the standard set
$C=\bigcup\limits_{i=1}^n(\o\xi_i,\o\eta_i)$. Then $dx(C)\mu_L(\st^{-1}(\*C))$. Obviously,\\ 
$\st^{-1}(\*C)\Delta
B(I)\sqe\bigcup\limits_{i=1}^n(\Mu(\o\xi_i)\cup \Mu(\o\eta_i))=\Mu(\partial C)$,
where the monad of a number $a\in[0,1]$ is denoted by $\Mu(a)$. Since the Loeb measure
of the monad of any number is equal to $0$ and so, $\Mu(\partial C)=0$, one has
$dx(C)=\o\left(\frac{|S(I)|}M\right)$. Substituting $[0,1]$ for $X$, $dx$ for $\nu$ and $\tau$ for $\psi$ in (\ref{hyp-iso3}) obtain
$dx(C)=dx(\tau^{-1}(C))=\mu_L(F^{-1}(\st^{-1}(C)))$. Since $I\setminus
F^{-1}(\st^{-1}(C))\sqe\Mu(\partial C)$, one has $\o\left(\frac{|I|}M\right)\leq
\o\left(\frac{|S(I)|}M\right)$. This means that if
$r_I=\max\{0,|I|-|S(I)|\}$, then $\frac{r_I}M\approx 0$. Let
$r=\max\{r_I\ | \ I\sqe Y \}$. Fix an arbitrary set
$Z\sqe\*[0,1]\setminus Y$ such that $|Z|=r$. For every $y\in Y$
set $S'(y)=S(y)\cup Z$ and for an arbitrary $I\sqe Y$ set
$S'(I)=\bigcup\limits_{y\in Y}S'(y)$. Then $S'(I)=S(I)\cup Z$,
$|S'(I)|=|S(I)|+r\geq |I|$, since $|I|-|S(I)|=r_I\leq r$. By the
Marriage Lemma there exists an injective map $\theta:Y\to
S'(Y)=Y\cup Z$ such that $\all\,y\ \theta(y)\in S'(y)$. Obviously
$|\theta^{-1}(Z)|=|Y\setminus\theta(Y)|\leq r$. So, there exists a
bijective map $\l:\theta^{-1}(Z)\to Y\setminus\theta(Y)$. Define
$T_{\d}:Y\to Y$ by the formula
$$
T_{\d}(y)=\left\{\begin{array}{ll} \theta(y),\ &y\in
Y\setminus\theta^{-1}(Z) \\ \l(y),\ &
y\in\theta^{-1}(Z)\end{array}\right.
$$
Notice that $\frac{|\theta^{-1}(Z)|)}M\leq \frac rM$. By
construction of $T_{\d}$ one has $\all\,y\in
Y\setminus\theta^{-1}(Z)\ |T_{\d}(y)-\tau(y)|<\d$.
Since $\mu_L(\theta^{-1}(Z))\leq\frac rM\approx 0$, the statement
(A) is proved.

Let $\cs(Y)$ be the set of all internal permutations of $Y$. Consider
the external function $f:\N\to \cs(Y)$ such that $f(n)=T_\frac
1n$. By $\aleph_1$-saturation the function $f$ can be extended to an internal function $\bar
f:\{0,\dots,K\}\to \cs(Y)$ for some $K\in\iN$. Internal function
$g(n)=\frac{|\{y\in Y\ |\ |F(y)-\bar f(y)|\geq\frac 1n\}|}M$ assumes
only infinitesimal values for all standard $n$. By Robinson's
Lemma there exists $L\in\iN$ such that $g(L)\approx 0$. set
$T=\bar f(L)$. Then $\mu_L(\{y\in Y\ |\ T(y)\approx F(y)\})=1$.
Since $F$ is a lifting of $\tau$, the same is true also for
$T(y)$. This proves I.

We have to prove now that a h.a. $T$ of $\tau$ can be chosen as
a cycle of maximal length.

II. Fix a permutation $T:Y\to Y$ that is a h.a. of $\tau$ and represent it
by a product of pairwise disjoint cycles, including the cycles of length 1 (fix points):
\begin{equation}\label{II1}
T=(y_{11}...y_{1n_1})(y_{21}...y_{2n_2})...(y_{b1}...y_{bn_b}),
\end{equation}
where $y_{ij}\in Y$ is the $j$-th element in the $i$-th cycle and $b$ is the number of cycles.
So,
\begin{equation}\label{II2}
\sum\limits_{i=1}^bn_i=M=|Y|.
\end{equation}
We assume also that $n_1\geq n_2\geq\dots\geq n_b.$
Consider the cycle
\begin{equation}\label{II3}
C=(y_{11}...y_{1n_1}y_{21}...y_{2n_2}...y_{b1}...y_{bn_b})
\end{equation}
By (\ref{II2}) $C$ is a cycle of length $M$, i.e. a transitive permutation.

Set $B=\{y\in Y\ |\ C(y)\neq T(y)\}$.
\begin{equation}\label{II4}
|B|=b=\sum\limits_{n=1}^Ma_n,
\end{equation}
where $a_n$ is the number of cycles of length $n$.

III.  Recall that a point $x\in [0,1]$ is said to be an $n$-periodic point of the transformation $\tau$
if its orbit under this transformation consists of $n$-points:
$x,\tau x,\dots,\tau^{n-1}x$. A point $x$ is said to be $\tau$-periodic if it is
$n$-periodic for some $n$.  The transformation $\tau$ is said to
be aperiodic if the set of periodic
points has measure zero. It is well-known that every measure
preserving automorphism $\tau$ of a Lebesgue space $X$ defines the
partition of this space by $\tau$-invariant Lebesgue subspaces of
aperiodic and $n$-periodic points. So, it is enough to prove our
statement for the case of aperiodic transformation $\tau$ and for
the case of $n$-periodic transformation $\tau$.

Suppose that the transformation $\tau$ is aperiodic. Let us
prove that under this assumption the cycle $C$ defined in
the part II is a h.a. of $\tau$.

Let $P_n(T)\sqe Y$ be the set of all $n$-periodic points of $T$ and let $P_n(\tau)\sqe X$ be the
set of all $n$-periodic points of $\tau$.
Since $T$ is a lifting of $\tau$
it is easy to that for every standard $k$ the following relations
\begin{equation}\label{IV1}
T(y)\approx\tau(\o y),\dots,T^{k}(y)\approx\tau^{k}{(\o y)}
\end{equation}
hold $\mu_L$-a.e. on $Y$.
So, for every standard $n$ $P_n(T)\sqe st^{-1}(P_n(\tau))$ up to a set of the Loeb measure zero. Since
$dx(P_n(\tau))=0$, one has $\frac 1M|P_n(T)|\approx 0$. Obviously, $|P_n(T)|=na_n$.
Thus, for every standard $n$ one has $\frac 1M\cdot
a_n\approx 0$.
By the Robinson's Lemma there exists an infinite $N$ such that $
\frac 1M\sum\limits_{n=1}^N a_n\approx 0.$. Obviously $
M\geq\sum\limits_{n=N+1}^Ma_n\cdot n\geq
(N+1)\sum\limits_{n=N+1}^M a_n.$ So, $\frac
1M\sum\limits_{n=N+1}^M a_n\leq\frac 1{N+1}\approx 0$ and $\frac
1M\cdot |B|=\frac 1M\sum\limits_{n=1}^Ma_n\approx 0. $ Thus,
$\mu_L(B)=0$, $C(y)=T(y)$ $\mu_L$-a.e. and $C$ approximates $\tau$.

IV. Suppose now that $\tau$ is $n$-periodic. We prove first that a
h.a. $T$ of $\tau$ also can be chosen to be $n$-periodic.
The relations (\ref{IV1}) imply that for almost every point $y\in Y$ if
$y$ has a standard period with respect to $T$, then this period is a multiple of $n$. Indeed, if
$y$ satisfies (\ref{IV1}), and its standard period is $nq+r$ for
$0<r<n$, then $\o y=\o T^{nq+r}(y)=\tau^{nq+r}(\o y)=\tau^{r}(\o
y)$, which is impossible since $\tau$ is $n$-periodic. By
$\aleph_1$-saturation, there exist an internal set $I\sqe Y$ such
that $\mu_L(I)=1$ and a number $N\in\iN$ such that for every point
$y\in I$, whose period is less, than $N$, this period is a
multiple of $n$.

Consider the representation (\ref{II1}) of $T$ and set $n_i=nq_i+r_i,\ r_i<n$ for each $i\leq b$.
Let $Y'\sqe Y$ be the set obtained by deleting from $Y$ the last $r_i$ elements of the i-th cycle
for each $i\leq b$.
The set $Y'$ has the Loeb measure equal to
$1$. Indeed, all the deleted elements either belong to the set
$Y\setminus I$, whose measure is $0$, or to a cycle whose length
is greater, than $N$. The number of these cycles does not exceed
$\frac MN$ and the number of deleted points in each such cycle is
less, than $n$. So the Loeb measure of the set of these points is
also equal to $0$.
Since $\mu_L(Y')=1$ the pair $(Y',\st)$ is a h.a. of $[0,1]$.
The construction of $Y'$ defines also the permutation $T':Y'\to Y'$
such that
\begin{equation}\label{V1}
T'=(y_{11}...y_{1\,n\cdot q_1})(y_{21}...y_{2\,n\cdot q_2})...(y_{b1}...y_{b\,n\cdot q_b}).
\end{equation}
Notice, that actually the number of cycles in $T'$ may be less, than $b$, since in case of $q_i=0$ the $i$-th cycle is empty.
However, the dynamical system $(Y',\mu_L,T')$ is a h.a. of the dynamical system $(X,\nu,\tau)$.
Indeed, let $D=\{y\in Y'\ |\ T(y)\neq T'(y)\}$. Then $D\sqe\{y\in Y'\ |\ T(y)\in Y\setminus Y'\}\sqe T^{-1}(Y\setminus Y')$.
Thus, $\mu_L(D)\leq \mu_L(Y\setminus Y')=0$
To obtain an $n$-periodic h.a. of $\tau$ it is enough to split each cycle in the representation (\ref{V1}) in cycles of length $n$.
Indeed, let the obtained cycle be
$$
T''=(z_1,...,z_n)(z_{n+1},...,z_{2\cdot n})\dots(z_{(i-1)\cdot n+1},...,z_{i\cdot n})\dots (z_{(K-1)\cdot n+1},...,z_{K\cdot n}),
$$
where $K=|Y'|/n$.
It is easy to see that $T''(y)\neq T'(y)$, only for the points
$z_{i\cdot n}$. Notice, that $\mu_L(\{z_{i\cdot n}\ |\ i\leq K\})=\frac 1n>0$
However, due to (\ref{IV1}) and the $n$-periodicity of $\tau$, for almost all of these points one has
$$T'(z_{i\cdot n})=z_{in+1}\approx\tau(\o z_{i\cdot n+1})=\tau^n(\o z_{(i-1)\cdot n+1})=\o z_{(i-1)\cdot n+1}.$$
At the same time $T''(z_{i\cdot n})=z_{(i-1)\cdot n+1}$ by the
definition. Thus, $T''(y)\approx T'(y)$ for almost all $y$.

V. To complete the proof of the theorem for $X=[0,1]$ we need to consider the case when all
orbits of $T$ have the same standard period $n$. In this case
$M=N\cdot n$.

It is easy to see that there exists a selector $I\subset Y$ (see the proof of Lemma \ref{0}) that is dense in
$\*[0,1]$, i.e. the monad $M(I)=\*[0,1]$. It is enough to show the
existence of a selector that intersects every interval with
rational endpoints. Obviously, for every finite set $A$ of such
intervals, there exists a selector that intersects each
interval from $A$. The existence of a dense selector follows from
the $\aleph_1$-saturation.

Let $I=\{y_1<y_2<\dots<y_N\}$ be a dense selector. Here $<$ is the
order in $\*[0,1]$. Due to the density of $I$ in $\*[0,1]$ for
every $k<N$ one has $y_k\approx y_{k+1}$. Obviously, the transformation $T$
can be represented by a product of pairwise disjoint cycles as follows:
$$
T=(y_1,...,T^{n-1}y_1)(y_2,...,T^{n-1}y_2)\dots(y_N,...,T^{n-1}y_N)
$$
Consider the following cycle $S$ of the length $M$:
$$
S=(y_1,...,T^{n-1}y_1y_2,...,T^{n-1}y_2\dots y_N,...,T^{n-1}y_N)
$$
Since for every $k\leq N$ holds $T^n(y_k)=y_k$, one has
$$\o S(T^{n-1}(y_k))=\o y_{k+1}=\o y_k =\o T^n(y_k)=\o
T(T^{n-1}(y_k))=\tau(\o T^{n-1}(y_k))$$ for almost all $k$. Thus,
$\o S(y)=\tau(\o y)$ for almost all $y$ and the cycle $S$ is a h.a. of $\tau$.

We proved actually that for every h.a. $Y$ of $([0,1]),dx)$ there exists an internal set $Y'\sqe Y$ with $\mu_L(Y')=1$
and a permutation $T':Y'\to Y'$ such that the hyperfinite dynamical system $(Y',\mu_L,T')$ is a h.a. of the dynamical system
$(X,\nu,\tau)$ and $T'$ is a transitive permutation of $Y'$ (see Part IV of this proof). To obtain a transitive h.a. $T:Y\to Y$ of $\tau$, set
$T'=(z_1,...,z_{|Y'|})$ and $Y\setminus Y'=\{u_1,...,u_{|Y\setminus Y'|}\}$ and consider the cycle of the length $|Y|$
$$
T:(z_1,...,z_{|Y'|},u_1,...,u_{|Y\setminus Y'|})
$$
Since $\mu_L(\{y\in Y\ |\ T'(y)\neq T(y)\})=0$ the transformation $T$ is h.a. of $\tau$.

VI. The statement of the theorem for the case of an arbitrary dynamical system $(X,\nu,\tau)$ satisfying the conditions of
Proposition \ref{Space-h-a} follows immediately from Lemma \ref{hyp-iso}.
Indeed, let a set $B\sqe X$, a set $C\sqe [0,1]$, a bijective map $\psi:B\to C$ and a bijective lifting $G:Y\to\*[0,1]$ of $\psi$
satisfy the conditions of Lemma \ref{hyp-iso}. Then $\lambda=\psi\tau\psi^{-1}:[0,1]\to [0,1]$ is a measure preserving
transformation. Fix an arbitrary h.a. $Y$ of the measure space $(X,\nu)$. Then by Lemma \ref{hyp-iso} the hyperfinite
set $Z=G(Y)$ is a h.a. of $([0,1],dx)$. By the results proved in the parts I-V, there exists a permutation $S:Z\to Z$
that is a h.a. of $\lambda$. Then it is easy to see that the permutation $T=G^{-1}SG:Y\to Y$ is a h.a. of $\tau$.
Obviously, if $S$ is a transitive permutation, then $T$ is a transitive permutation as well. $\Box$

The following corollary of Theorem \ref{cycle-appr} is obvious.

\begin{Cor} \label{kam-cor} For every Lebesgue dynamical system $(X,T,\nu)$, there exists a transitive
Loeb dynamical system $(Y, T)$ and a measure preserving map $\psi:Y\to X$.\end{Cor}

In other words, every Lebesgue dynamical system is a homomorphic image of an appropriate transitive Loeb dynamical
system.

\textbf{Proof}. It is enough to take any transitive h.a. of $(X,T,\nu)$ for $(Y,T)$ and $st:Y\to X$ for the map $\psi$ $\Box$

Corollary \ref{kam-cor} immediately allows to deduce the ergodic theorem for Lebesgue spaces from
ergodic theorem for Loeb space. This was done in the paper \cite{Kam}. Theorem 2 of \cite{Kam} is stronger
than Corollary \ref{kam-cor}. It states that each Lebesgue dynamical is a homomorphic image of \textbf{any} transitive Loeb dynamical
system.

On the other hand Theorem 2 of \cite{Kam} does not imply (at least, immediately) the existence of a hyperfinite
approximation of a Lebesgue dynamical system. The notion of a hyperfinite approximation has a simple standard interpretation
(see the end of this section) and might be of interest by itself, in particular, from the point of view of computer simulations of dynamical systems.
Hyperfinite approximations of many important dynamical systems can be easily described (see examples 5 -- 7 below), while Loeb preimages of these systems
obtained by Theorem 2 of \cite{Kam} do not have any simple description.

In what follows we say that a number $N\in\*\N$ is \emph{$M$-bounded}, if
the number $\frac NM$ is bounded. For any $y\in Y$ and any $k\in\*\N$ the set $\{T^ky,T^{k+1}y,...,T^{k+N-1}y\}$ is said
to be an $N$-segment of the orbit of $y$. The $N$ segment $\{y,Ty,...,T^{N-1}y\}$ is said to be initial.

It is easy to see that even a transitive dynamical system $(Y,T)$ is never
ergodic. Indeed, for any $y\in Y$ and $N\in\iN$ such that
$0<\o\left(\frac NM\right)<1$ every $N$-segment $A$ of the orbit of
$y$ is obviously a $T$-invariant set, since $\mu_L(A\bigtriangleup T(A))=0)$.
Every transitive system is equivalent to a system considered in the examples 1-4. Examples 3
and 4 show, that for every transitive system $(Y,\mu_L, T)$ there
exist an $S$-integrable function $F$ and numbers $N,K\in\iN$ such
that $\frac NM,\frac KM$ are bounded and $A_N(F,T,y)\not\approx
A_K(F,T,y)$ for all $y$ in some set of a positive Loeb measure.

In the few following propositions and examples 5 - 7 we discuss the behavior of ergodic means $A_N(F,T,y)$ for $M$-bounded $N\in\iN$ in
the case, when a Loeb dynamical system $(Y,\mu_L,T)$ is a h.a. of a Lebesgue dynamical system $(X,\nu,\tau)$ and $F$ is an $S$-integrable
lifting of a function $f\in L_1(\nu)$.

The following proposition is an easy corollary of Theorem \ref{NSBET}.

\begin{Prop} \label{CNSBET} In conditions of the previous paragraph let $\tilde{f}=\lim\limits_{n\to\infty}A_n(f,\tau, x)$ and
let $\widetilde{F}$ be an $S$-integrable lifting of $\tilde{f}$, then there exists an
$N\in\iN$ such that $\mu_L$-a.e. $\all\, K\in\iN\ (K<N\Rightarrow A_K(F,T,x)\approx\widetilde{F}(x))$.
\end{Prop}

\begin{Cor} \label{non-erg} Let $T$ be a transitive permutation and let $\tau$
be a non-ergodic transformation. Consider a function
$f\in L_1(\nu)$ such that the set $B\sqe X$ of all $x\in X$ satisfying
inequality $\lim\limits_{n\to\infty}A_n(f,\tau,x)\neq Av(f)$ has a
positive measure $\nu$. Then there exist infinite
$M$-bounded $N,K$ such that for almost all $y\in\st^{-1}(B)$ one
has $A_N(F,T,y)\not\approx A_K(F,T,y)$.
\end{Cor}

\textbf{Proof} \newcommand{\wt}{\widetilde}. Let $\wt
f=\lim\limits_{n\to\infty}A_n(f, \tau,\cdot)$ and $\wt F$ be the same
as in Proposition \ref{CNSBET}. By this proposition there exists $N\in\iN$
such that $\frac{N}{M}\approx 0$ and $A_N(F,T,y)\approx\widetilde{F}(y)$ $\mu_L$-a.e.
Thus, $A_N(F,T,y)\not\approx Av(f)$ for $\mu_L$-almost all $y\in\st^{-1}(B)$.

On the other hand, since $T$ is a cycle of length $M$, by Theorem \ref{ErgMeanStab} one has
$A_K(F,T,y)\approx Av(F)\approx Av(f)$ for all $y\in Y$ and for all $K$ such that $\frac{K}{M}\approx 1$. Thus,
$A_K(F,T,y)\not\approx A_M(F,T,y)$ for $\mu_L$-almost all $y\in\st^{-1}(B)$.
$\Box$

Let $X$ be a compact metric space. Consider a hyperfinite set $Y\sqe \*X$. This set defines a
Borel measure $\nu_Y$ on $X$ by the formula $\nu_Y(K)=\mu_L(st^{-1}(K\cap Y))$. Obviously $Y$ is a h.a.
of the measure space $(X,\nu_Y)$. Let $T:Y\to Y$ be an internal permutation that is
$S$-continuous on $A$ 
for some (not necessary internal) set
$A\sqe Y$ with $\mu_L(A)=1$, i.e.
\begin{equation} \label{measurable}
\all\, a_1,a_2\in A\ (a_1\approx a_2\Limpl T(a_1)\approx T(a_2)).
\end{equation}
Notice that since $\st^{-1}(st(A))\supseteq A$ and $\mu_L(A)=1$,
the set $\st(A)\sqe X$ is a measurable set w.r.t. the completion of the measure $\nu_Y$, which we denote by $\nu_Y$ also,
and $\nu_Y(\st(A))=1$.

Define a map $\tau_T:X\to X$ such that $\tau_T(\st(y))=\st(T(y))$ for $y\in A$ and $\tau_T\upharpoonright{X\setminus\st(A)}$
is an arbitrary measurable permutation of the set $X\setminus\st(A)$.
\begin{Prop} \label{tauT}
The map $\tau_T$ preserves the measure $\nu_Y$.
\end{Prop}
\textbf{Proof}. Replacing, if necessary, $A$ by $\bigcap\limits_{n\in\N}T^n(A)$ we may assume that $A$ is invariant for permutation $T$.
Then, obviously, $\st(A)$ is invariant for $\tau_T$.

Consider a closed set $B\sqe X$. We have to prove that $\nu_Y(\tau^{-1}_T(B))=\nu_Y(B)$. One has
$$\nu_Y(\tau_T^{-1}(B))=\nu_Y(\tau^{-1}_T(B)\cap\st(A)).$$
It is easy to check that
$$\tau_T^{-1}(B)\cap\st(A)=\st(T^{-1}(\st^{-1}(B))\cap A).$$
Thus,
$$\nu_Y(\tau_T^{-1}(B))=\mu_L\left(\st^{-1}(\st(T^{-1}(\st^{-1}(B))\cap A))\right)=\mu_L\left(\left(\st^{-1}(\st(T^{-1}(\st^{-1}(B))\cap A))\right)\cap A\right)$$
Using (\ref{measurable}) and the $T$-invariance of it is easy to check, that
$$\st^{-1}(\st(T^{-1}(\st^{-1}(B))\cap A))\cap A=T^{-1}(\st^{-1}(B))\cap A).$$
So,
$$\nu_Y(\tau_T^{-1}(B))=\mu_L(T^{-1}(\st^{-1}(B))\cap A))=\mu_L(T^{-1}(\st^{-1}(B))=\mu_L(\st^{-1}(B))=\nu_Y(B). $$
In the last chain of equalities we used the facts that $\mu_L(A)=1$ and that $T$ being a permutation preserves the Loeb measure. $\Box$

\begin{Prop}\label{functional}
1) In conditions of Proposition \ref{tauT} for any $a>0$ and for any $y\in Y$ the following positive functional $l_a(\cdot,T,y)$ on $C(X)$ is defined:
$l_a(f,T,y)=\o A_K(\*f,T,y)$, where $\o\left(\frac KM\right)=a$ and $f\in C(X)$

2) If $\all\,K,L\in\iN\ (\frac KM\approx\frac LM\approx 0\Limpl A_K(\*f,T,y)\approx A_L(\*f,T,y)$), then $l_0(f,T,y)$ is defined by the same formula as in 1).
In this case $l_0(f,T,y)=\wt f(y)$

3) If $T:Y\to Y$ is $S$-continuous, then the functional $l_0(\cdot,T,y)$ is $\tau_T$-invariant for all $y\in Y$.
\end{Prop}

\textbf{Proof}. The correctness of the definition in 1) follows from Theorem \ref{ErgMeanStab}. The statement 2) follows from Proposition \ref{WeakNSBET}. To prove
statement 3) notice that if $T$ is $S$-continuous on $Y$, then $\tau_T$ is continuous on $X$ and, thus, $\*(f\circ\tau_T)\upharpoonright Y$ is a lifting of $f\circ\tau_T$. So,
$$\all\, y\ \all\, K\in\*\N\ \*(f\circ\tau_T) (T^K(y))\approx f(\tau_T(\o T^K(y)))=f(\o T^{K+1}(y))\approx\* f(T^{K+1}y)$$
These equivalences allows to prove that $A_K(\*(f\circ\tau),T,y)\approx A_K(\* f,T,y)$ (see the proof of Proposition \ref{S-int-means}). $\Box$

Recall that a continuous transformation $\tau:X\to X$ is said to
be \emph{uniquely ergodic} if there exists only one
$\tau$-invariant Borel measure on $X$\footnote{Krylov-Bogoljubov
theorem states the existence of at least one $\tau$-invariant
measure.}.

\begin{Th} \label{un-erg} If $\tau$ is a uniquely ergodic transformation of a compact metric space $X$,
$Y\sqe \*X$ is a hyperfinite set such that $\st(Y)=X$, and $T:Y\to Y$ is an internal permutation such that $\all\,y\in Y\ \st (T(y))=\tau(\st (y))$,
then for every $y\in Y$ such that the $\tau$-orbit of $\st(y)$ is dense in $X$, for every $N\in\iN$ and for every $f\in C(X)$ one has
\begin{equation} \label{un-erg-1}
A_N(\*f\upharpoonright Y,T,y)\approx\int\limits_Xfd\nu,
\end{equation}
where $\nu$ is the $\tau$-invariant measure.
\end{Th}

\textbf{Proof}. Let $y\in Y$ satisfy conditions of the theorem. For a number $K\in\*\N$ 
we denote the initial $K$-segment of the $T$-orbit of $Y$ by $S(K,y)$. Then
for any $K\in\iN$ one has $\st(S(K,y))=Y$, since the closed set $\st(S(K,y))$ contains the $\tau$-orbit of $\st(y)$. Let $K$ be a $T$-period of $y$. Then $K\in\iN$. Otherwise, the $\tau$-orbit of $\st(y)$ would be finite, while we assume $X$ to be infinite. It is easy to see that it is enough to prove the theorem
for every $N\in\iN$ such that $N\leq K$. Under this assumption all elements of the set $Y_1=\{y, Ty, ... , T^{N-1}y\}$ are distinct. Since $\st(Y_1)=X$, the set $Y_1$
defines the Borel measure $\nu_{Y_1}$ on $X$. Let $T_1:Y_1\to Y_1$ be the permutation of $Y_1$ that differs from $T$ only for one element $T^{N-1}y$: $T_1(T^{N-1}y)=y$.
Set $A=Y_1\setminus\{T^{N-1}y\}$. Then $X$, $\tau$, $Y_1$, $T_1$, and $A$ satisfy conditions of Proposition \ref{tauT}: $\mu_L(A)=1$,
$\all\, z\in A\ \st(T_1z) =\tau(\st(z))$, i.e. $\tau_{T_1}=\tau$ and $T_1$ is $S$-continuous on $A$, since $\tau$ is a continuous map. By Proposition \ref{tauT}
the measure $\nu_{Y_1}$ is $\tau$-invariant. Thus, $\nu_{Y_1}=\nu$ due to the unique ergodicity of the map $\tau$. If $f\in C(X)$, then obviously $\*f\upharpoonright Y_1$ is an $S$-integrable lifting of $f$. This proves the equality (\ref{un-erg-1}). $\Box$

\textbf{Example 4}. Let $Y=\{0,\frac 1m,\frac 2M,\dots,\frac
{M-1}M\}$ and $T:Y\to Y$ is defined by the formula $T\left(\frac
kM\right)=\frac{k\oplus 1}M$, where $\oplus$ is the addition modulo
$M$. Then $(Y,\st)$ is a hyperfinite approximation of the segment
$[0,1]$ with the Lebesgue measure, and the map $T$ is a h.a. of
the identical map $id:[0,1]\to [0,1]$. Obviously, the dynamical
system $(Y,\mu_L,T)$ is isomorphic to the one, considered in
Examples 1 - 3 up to the trivial isomorphism $k\mapsto\frac kM$.
The function $F(k)$ of Example 3 is a lifting of the function
$f(x)=x$, which is the function $f_0$ of Example 3. Obviously
$\all\,n,x\ A_n(f,id,dx)=f(x)$. So, in this example
$\all\, N\in\iN\ (\frac{N}{M}\approx 0\Rightarrow
A_N(F,T,y)\approx \widetilde{F}(y))$ $\mu_L$-a.e.

\textbf{Example 5}. Let $\a\in\R$. Denote by $S(\a)$ the shift of the interval
$[0,1]$ by $\a$, i.e. the dynamical system $D=([0,1],\tau,dx)$, where $\tau(x)=x\oplus\a$
and $\oplus$ is the addition modulo $1$.
In this example we construct two distinct h. a. of the dynamical system $S(\frac 12)$. Similar considerations are appropriate for an
arbitrary rational shift of the unit interval.

Let $Y$ be the same as in the previous example and $M=2N$ and let $T:Y\to Y$ be given by the formula: $T(y)=y\oplus\frac NM$. It is easy to see
that  $\all\,y\in Y\ T(y)=\tau(y)$ and if $y\not\approx \frac 12$, then $\tau(\o y)=\o\tau(y)$. Since $\mu_L(Y\cap\Mu(\frac 12))=0$,
the dynamical system $(Y,\mu_L,T)$ is a h.a. of the system $S(\frac 12)$. It is easy to see that for an arbitrary $S$-integrable function $F$ on $Y$, for
every $N\in\iN$ and for every $y\in Y$ one has $A_N(F,T,y)\approx\frac 12 (F(y)+F(y\oplus\frac 12))$ (=, if $N$ is an even number).
In particular, if $F$ is an $S$-integrable lifting of a function $f\in L_1(0,1)$ then for almost all $y\in Y$  and for all $N\in\iN$ one has
\begin{equation} \label{11}
A_N(F,T,y)\approx\frac 12(f(\o y)+f(\o y+\frac 12))=\lim\limits_{n\to\infty}A_n(f,\tau,\o y)=\wt f(\o y).
\end{equation}

In this case every point $y\in Y$ is periodic with the period $2$. By Theorem \ref{cycle-appr} there exists a transitive h.a. of $S(\frac 12)$. In this case
it can be easily constructed explicitly.

Let $(Y,\mu_L,T)$ be the same as above, but $M=2N+1$. Then it is easy to see that if $y=\frac kM$, then $T(y)=\frac{k+N\,(\mod\,M)}{M}$. So, $T$ is isomorphic to a permutation of $\{0,1,...,M-1\}$ given by the formula $k\mapsto k+N\,(\mod\,M)$. Since $\gcd(N,M)=1$, this permutation is transitive.
It is easy to see that the transformation $T$ is a h.a. of $\tau$. Indeed, consider the set $B=[0,1)\setminus \{\frac 12\}$. Then $\tau$ is continuous on
$B$. Let $A=st^{-1}(B)=Y\setminus(\Mu(0)\cup\Mu(1)\cup\Mu(\frac 12))$. Then $\mu_L(A)=1$ and $A$ is invariant for $T$. An easy calculation show
that $\all\,y\in A\ \*\tau(y)-T(y)=\frac 1{2M}$. Since $\*[0,1]\setminus(\Mu(0)\cup\Mu(1)\cup\Mu(\frac 12))$ is invariant for $*\tau$,  one has
$T(y)\approx\tau(y)\approx\tau(\o y)$. Moreover, it is easy to see that, if $y=\frac kM\in A$, then for every $K\in\*\N$ such that $\frac{k-K}M\in A$ one has

\begin{equation} \label{12}
\*\tau^K(y)-T^K(y)=\frac K{2M}
\end{equation}

So, if $\frac KM\approx 0$, then $\all\, n\leq K$ the set $A$ is invariant for $T^n$ and $\all\, y\in Y\ \o T^n(y)\approx\tau^n(\o y)$ $\mu_L$-a.e.

Let $f:[0,1]\to\R$ be a standard continuous function, then obviously $\*f\upharpoonright Y$ is an $S$-integrable lifting of
$f$. Assume first that $f$ satisfies the Lipschitz condition. Then, using the equality (\ref{12}), one immediately obtains that if $\frac KM\approx 0$, then
$\all\, y\in A\ A_K(\*f,T,y)\approx A_K(\*f,\tau,y)\approx \wt f(\o y)$ (see formula (\ref{11})).
For any $K\in\*\N$ denote $l_K(\cdot,T,y)$ the positive functional on $C(X)$ given by the formula $l_K(f,T,y)=\o A_K(\*f, T,y)$.
Since in the case of $\frac KM\approx 0$ for any function $f$ satisfying the Lipschitz condition $l_K(f,T,y)=\wt f(\o y)$ and such functions are dense in $C(X)$,
we see that the functional $l_0(\cdot, T, y)$ of Proposition \ref{functional} is defined. The formula (\ref{11}) shows that
$$l_0(\cdot, T,y)=\frac 12\left(\d_{\o y}+\d_{\o y\oplus 0.5}\right)$$

Using formula (\ref{12}) one can easily obtain explicit expressions for $l_a(\cdot, T, y)$ for positive $a\in\R$. These expressions depend on relations
between $a$ and $\o y$. For example, if $0.5a<\o y<0.5$, then
$$
l_a(f,T,y)=\frac 1a \left[\int\limits_{\o y-0.5a}^{\o y}f(t)dt+\int\limits_{\o y+0.5-0.5a}^{\o y+0.5}f(t)dt\right].
$$
We see that $l_a\to \l_0$ as $a\to 0$. It easy to see also that for all $a>0$ and $y\in A$ the functional $l_a$ is $\tau$-invariant, though $\tau$ is not continuous
everywhere on $[0,1]$. Compare with Proposition \ref{functional}

{\bf Example 6}. Let $\tau$ be a shift of $[0,1]$ by an irrational number $\a>0$ modulo 1. It is well known
that $\tau$ is uniquely ergodic.  Consider $M,N\in\iN$ such that
$\gcd(M,N)=1$ and $\frac NM\approx\a$. Consider the same finite space $Y$ as in Example 5.
Let $T:Y\to Y$ be given by the formula $T(y)=y+\frac NM\,(\mod 1)$. By the same reasons as in
Example 6 the permutation $T$ is a cycle of length $M$ and $T$ is a h.a. of $\tau$. Thus, for every $S$-continuous function
$F:Y\to\*\R$, $K\in\iN$ and $y\in Y$ one has $A_K(F,T,Y)\approx Av(F)$. Compare with Example 5.

{\bf Example 7}. (Approximations of Bernoulli shifts). Let $\Sigma_m=\{0,1,\dots m-1\}$. Consider the compact space $X=\Sigma_m^{\mathbb{Z}}$ with the Tychonoff topology. Let $a$ be a function, such that $\dom(a)\subset\Z$ is finite, and $\range(a)\sqe\Sigma_m$. Let $S_a=\{f\in X\ |\ f\upharpoonright\dom(a)=a\}$. Then the family of all such $S_a$ form a base of neighborhoods of the compact space $X$. For $g\in\*X$ set $f=g\upharpoonright\Z$, then $f\in X$ and it is easy to see that $f=\st(g)$.

The continuous transformation $\tau:X\to X$ defined by the formula $\tau(f)(n)=f(n+1)$ where $f\in X$ and $n\in\Z$ is an invertible Bernoulli shift.
Every probability distribution $\{p_0,\dots,p_{m-1}\}$ ($p_i>0,\ \sum\limits_{i=0}^{m-1}p_i=1$) on $\Sigma_m$ defines a Borel measure on $X$ that is obviously
invariant with respect to $\tau$. It is well-known that $\tau$ is ergodic for each of these measures. So, the transformation $\tau$ is not uniquely ergodic.
Here we restrict ourselves only to the case of the uniform distribution on $\Sigma_m$, i.e. to the case of $p_0=\dots=p_{m-1}=\frac 1m$.
The corresponding Borel measure on $X$ is denoted by $\nu$. Obviously, $\nu(S_a)=m^{-|\dom(a)|}$.

We construct here two hyperfinite approximations of the dynamical system $(X,\nu,\tau)$. First we consider the straightforward approximation by a hyperfinite shift.
Fix $N\in\iN$ and set $Y=\Sigma_m^{\{-N,\dots,N\}}$. Then $M=|Y|=m^{2N+1}$. Define $\l:Y\to\*X$ as follows. For $y\in\Sigma_m^{\{-N,\dots,N\}}$ set
\begin{equation} \label{3}
\l(y)(n)=\left\{\begin{array}{ll} y(n),\ &|n|\leq N \\ 0,\ &|n|>N\end{array}\right.
\end{equation}
Then $\st\circ\l (y)=y\upharpoonright\Z$. Thus, for every standard neighborhood $S_a$ 
defined above one has $(st\circ\l)^{-1}(\*S_a)=\{y\in Y\ |\ y\upharpoonright\dom(a)=a\}$. 
So, $\mu_L((\st\circ\l)^{-1}(\*S_a))=\nu(S_a)=m^{-\dom(a)}$. 
This proves that $(Y,\l)$ is a h.a. of $(X,\nu)$.

Certainly, an arbitrary internal map from $Y$ to $\Sigma_m^{\*\Z\setminus\{-N,\dots,N\}}$
can be used to define the values $\l(y)(n)$ for $|n|>N$ and $y\in Y$ in the definition of $\l$ (\ref{3})

In what follows we use notations $y_1\approx y_2$ and $\st(y)$ for $\l(y_1)\approx\l(y_2)$ and $\st(\l(y))$ respectively.

Define the map $S\to S$ by the formula $S(y)(n)=y\left(n+1(\mod 2N+1)\right)$ for any $y\in Y$ and $n\in\{-N,\dots,N\}$. Then obviously $\tau(\st(y))=\st(S(y))$ for all $y\in Y$. So, $(Y,\l,S)$ is a h.a. of the dynamical system $(X,\nu,\tau)$.

Since every point $y\in Y$ is $(2N+1)$-periodic with respect to $S$ the permutation $S$ is not transitive. Though the existence of a transitive h.a. of $\tau$ is proved in Theorem \ref{cycle-appr}, it is not easy to construct such an approximation explicitly.

To do this we reproduce here the construction of de Bruijn sequences.

\begin{Def} \label{de Bruijn}
An $(m,n)$-de Bruijn sequence on the alphabet $\Sigma_m$ is a sequence $s=(s_0,s_1,\dots,s_{L-1})$ of $L=m^n$ elements $s_i\in\Sigma_m$ such
that all consecutive subsequences $(s_i, s_{i\oplus 1},\dots,s_{i\oplus n-1})$ of length $n$ are distinct.

Here and below in this example the symbols $\oplus$  and $\ominus$ denote + and -  modulo $L$, so that the sequence $s$ is considered as a sequence of symbols from
$\Sigma_n$ placed on a circle.
\end{Def}

It was proved \cite{deB1, deB2} that there exist $(m!)^{m^{n-1}}\cdot m^{-n}$ $(m,n)$-de Bruijn sequences. See also \cite{deB3} for a simple algorithm for de Bruijn sequences and more recent references.

To construct a transitive h.a. $T:Y\to Y$ of $\tau$ fix arbitrary $(m, 2N+1)$ de Bruijn sequence $s=(s_0,s_1,\dots,s_{M-1})$ here $L=M$.
Let $y=(y_{-N},\dots,y_{-N})\in Y$. Then there exists the unique consecutive subsequence $\s(y)=(s_i,s_{i\oplus 1},\dots,s_{i\oplus 2N})$ such that
$y_j=s_{j\oplus i\oplus N})$.  Set $P(\s(y))=(s_{i\oplus 1},\dots, s_{i\oplus 2N\oplus 1})$ and $T(y)=\s^{-1}(P(\s(y)))$. Notice that if $i< M-N$ , then
for all $j\leq N$ one has $s_{j\oplus i\oplus N}=s_{j+i+N}$. So, $T(y)_j=y_{j+1}$ for all $j<N$ and, thus, for all standard $j$.
This last equality implies that $\st(T(y))=\tau(\st(y))$ for all $y\in Y$ such that the first entry of the sequence $\s(y)$ is the $i$-s term of the initial de Bruijn
sequence for $i \leq L-2N-1$. So, $\mu_L(\{y\ |\ \st(T(y))=\tau(\st(y))\})\geq\frac{M-N}M\approx 1$. This proves that $T$ is a h.a. of $\tau$. We call $T$ a de Bruijn approximation of $\tau$.

It is interesting to study the behavior of ergodic means of described approximations. This problem will be discussed in another paper. We confine ourselves with two simple remarks.

1. If $\s(y)=\la s_i,\dots,s_{i+2N}$ and $i<M-N$, then $A_n(F,T,y)=A_n(F,S,y)$ for all $n<N$.

2. Let $S_0=\{f\in X\ |\ f(0)=1\}$, so that $\nu(S_0)=\frac 12$ and let $\chi_0$ be a characteristic function of $S_0$. For $y\in Y$ let $f=\st(y)$. Set $A(y)=f^{-1}(\{1\})\cap\N)$. Recall that the density of $A(y)$ is given by the formula
$$d(A(y))=\lim\limits_{m\to\infty}\frac{|A(y)\cap\{0,\dots, m-1\}|}m$$
It is easy to see that for $m<N$ one has
$$A_m(\*\chi_0,T,y)=\frac{|A(y)\cap\{0,\dots, m-1\}|}m.$$
So, for all $y\in Y$ such that the density $d(A(y))$ exists one has $\ex\, K\in\iN\,\all\, m\in\iN\ (m\leq K\Limpl A_m(\chi_0,T,y)\approx d(A(y))$.
Due to Proposition \ref{CNSBET} there exist $K\in\iN$ such that for $\mu_L$-almost all $y\in Y$ one has $A_m(\*\chi_0, T,y)\approx\frac 12$.

\bigskip

\emph{Standard version of the notion of a hyperfinite approximation of a dynamical system.}

\bigskip

We use the same approach as above to formulate a sequence version of the notion of a hyperfinite approximation of a dynamical system.

\begin{Def} \label{seqap} Let $(X,\rho)$ be a compact metric space, $\nu$ be a Borel measure on 
$X$, $\tau:X\to X$ be a measure preserving transformation of $X$, $\{Y_n\sqe X\ |\ n\in\N\}$ 
be a sequence of finite approximation of the measure space $(X,\nu)$ in the sense of 
Definition \ref{st-hypap} and $T_n:Y_n\to Y_n$ be a sequence of permutations of $Y_n$.
We say that a sequence $T_n$ is an approximating sequence of the transformation $\tau$ if for every $N\in\iN$ the internal permutation $T_N:Y_N\to Y_N$ is a h.a. of
$\tau$ in the sense of Definition \ref{hypapds}. In this case we say that the sequence of finite dynamical systems $(Y_n,\mu_n, T_n)$ approximates the dynamical system
$(X,\nu,\tau)$. Here $\mu_n$ is a uniform probability measure on $Y_n$. \end{Def}

The reformulation of this definition in full generality in standard mathematical terms is practically unreadable. However, it is easy to reformulate it
for the case of an almost everywhere continuous transformation $\tau$. This case covers a lot of important applications.

Denote the set of all points of continuity of the map $\tau:X\to X$ by $D_{\tau}$.

\begin{Lm} \label{seqap1} Suppose that $\nu(D_{\tau})=1$ and let $Y\sqe X$ be a h.a. of the measure space $(X,\nu)$. Then a permutation $T:Y\to Y$ is a h.a. of the
transformation $\tau$ if and only if for every positive $\e\in\R$ one has
\begin{equation} \label{seqap2}
\frac 1{|Y|}\left({|\{y\in Y\ |\ \rho(T(y),\*\tau(y))>\e\}|}\right)\approx 0.
\end{equation}\end{Lm}

\textbf{Proof} ($\Limpl$) Let $A=\{y\in Y\ |\ \o y\in D_{\tau}\}=\st^{-1}(D_{\tau})$, $B=\{y\in Y\ |\ T(y)\approx\tau(\o y)\}$. Then, $\mu_L(A)=1$ since $Y$ is a h.a. of
the measure space $(X,\nu)$ and $\nu(D_{\tau})=1$. Since $T$ is a h.a. of $\tau$, one has $\mu_L(B)=1$. Thus, $\mu_L(A\cap B)=1$. Since $\tau$ is continuous on
$\D_{\tau}$ one has
\begin{equation}\label{seqap3}
\all\,x\in X\ \o x\in D_{\tau}\Limpl \*\tau(x)\approx\tau(\o x).
\end{equation}
So, $\all\, y\in A\cap B\ \*\tau(y)\approx\tau(\o y)$ and thus,  $\all\, y\in A\cap B\ \*\tau(y)\approx T(y)$. So, for every positive $\e\in\R$ one has
$\{y\in Y\ |\ \rho(T(y),\*\tau(y))>\e\}\sqe Y\setminus(A\cap B)$. This proves (\ref{seqap2}).

$(\Longleftarrow)$ Suppose that (\ref{seqap2}) holds for every positive $\e\in\R$. Then obviously $\mu_L(\{y\in Y\ |T(y)\approx\*\tau(y)\})=1$. On the other hand,
by (\ref{seqap3}) one has $\mu_L(\{y\in Y\ |\ \*\tau(y)\approx\tau(\o y)\})=1$. Thus, $\mu_L(\{y\in Y\ |\ T(y)\approx\tau(\o y)\})=1$, i.e. $T$ is a h.a. of $\tau$
$\Box$

Lemma \ref{seqap1} implies immediately the following

\begin{Prop}[Standard version of Definition \ref{seqap}] \label{seqap4} In conditions of Definition \ref{seqap} and Lemma \ref{seqap1} the sequence of permutations
$T_n:Y_n\to Y_n$ is an approximating sequence of the transformation $\tau$ if and only if for every positive $\e\in\R$ one has
\begin{equation}\label{seqap5}
\lim\limits_{n\to\infty}\frac 1{|Y_n|}\left({|\{y\in Y_n\ |\ \rho(T_n(y),\tau(y))>\e\}|}\right)=0.\qquad \Box
\end{equation}
\end{Prop}

\begin{Th}\label{exseqap} Let $(X,\rho)$ be a compact metric space and $\nu$ be a Borel measure on $X$ such that the measure space $(X,\nu)$ satisfies the conditions
of Proposition \ref{Space-h-a}. Then for every measure preserving transformation $\tau:X\to X$ with $\nu(D_{\tau})=1$ there exist a sequence of finite sets
$Y_n\sqe X$ and a sequence of permutations $T_n:Y_n\to Y_n$ such that the sequence of finite dynamical systems $(Y_n,T_n)$ approximates the dynamical system $(X,\nu,\tau)$
in the sense of Definition \ref{seqap}. Moreover, one can choose transitive permutations $T_n$.
\end{Th}

\textbf{Proof.} Let $Y_n\sqe X$ be a sequence that approximates the measure space $(X,\nu)$ in the sense of Definition \ref{st-hypap}. Such sequence exists by
Proposition \ref{ex-st-hypap}. Then by Proposition \ref{equivalence} for any $N\in\iN$ the set $Y_N$ is a h.a. of the measure space $(X,\nu)$ in the sense of Definition
\ref{hypap}. By Theorem \ref{cycle-appr} there exists a (transitive) permutation $T_N:Y_N\to Y_N$ that is a h.a. of the transformation $\tau$.
By Lemma \ref{seqap1}, since $\nu(D_{\tau})=1$, this means that $(Y_N,T_N)$ satisfies (\ref{seqap2}) for every standard positive $\e$. In this proof the
letter $T$ maybe with lower indexes always denotes a (transitive) permutation.

For every numbers $n,m\in\N$ define the set
$$A_{n,m}=\left\{k\in\N\ \left| \right.\ \ex\, T:Y_k\to Y_k\
\left(\frac 1{|Y_k|}\cdot\left| \{y\in Y_k\ |\ \rho(T(y),\tau(y))>\frac 1n\}\right|<\frac 1m\right)\right\}.$$

Since $\all\,N\in\iN\ N\in\*A_{n,m}$, there exists a standard function $N(n,m)$ such that $\all\,k>N(n,m)\ k\in A_{n,m}$. By the definition of sets $A_{m,n}$, there exists
a standard function $T(k,n,m):Y_k\to Y_k$ with the domain $\{(n,m,k)\in\N^3\ |\ k>N(n,m)\}$ such that
$$
\frac 1{|Y_k|}\cdot\left| \left\{y\in Y_k\ |\ \rho(T_k(y),\tau(y))>\frac 1n\right\}\right|<\frac 1m.
$$

Now it is easy to see that if $r=N(n,n)+n$, then the sequence $(Y_r,T_r)$ satisfies the conditions of Proposition \ref{seqap4} $\Box$

\bigskip

IICO-UASLP

AvKarakorum 1470

Lomas 4ta Session

SanLuis Potosi SLP 7820 Mexico

Phone: 52-444-825-0892 (ext. 120)

e-mail:glebsky@cactus.iico.uaslp.mx

\bigskip

Department of Mathematics and Computer Science

Eastern Illinois University

600 Lincoln Avenue

Charleston, IL 61920-3099 USA

Phone: 1-217-581-6282

e-mail: yigordon@eiu.edu

\bigskip

Department of Mathematics

University of Illinois at Urbana-Champaign

1409 W. Green Street

Urbana, Illinois 61801-2975 USA

e-mail: henson@math.uiuc.edu

\end{document}